\documentclass[12pt]{amsart}
\usepackage{amssymb}
\usepackage{eucal}
\usepackage{amsmath}
\usepackage{plesty}
\usepackage{paralist}
\usepackage[all]{xy}
\usepackage{enumitem}
\usepackage{tikz}

\newcommand{\Title}{Regular Decomposition of  Ordinarity in Generic Exponential Sums}

\begin{document}

\title{\Title}
\date{August 2012}

\author{Phong Le}
\address{Department of Mathematics, Duneavy Hall Rm \#333, P.O. Box 2044, Niagara University, New York}
\email{ple@niagara.edu}
\maketitle
\begin{abstract}

In \cite{WAN93} and \cite{WAN04} Wan establishes a decomposition theory for the generic Newton polygon associated to a family of $L$-functions of $n$-dimensional exponential sums over finite fields.  In this work we generalize the star, parallel hyperplane and collapsing decomposition, demonstrating that each is a generalization of a complete regular decomposition.  


\end{abstract}


\section{Introduction}






In \cite{AS89} Adolphson and Sperber established a combinatorial lower bound for the generic Newton polygon attached to a family of $L$-functions of $n$-dimensional exponential sums over finite fields.  Notably this lower bound is independent of the character of the finite field.  In the case of toric hypersurfaces, the bound required the use of Hodge numbers.  For this reason, the bound was called the Hodge Polygon.  Generic Newton polygons which coincide with the Hodge polygon are called generically ordinary.  In \cite{AS89} they also conjectured conditions when generic ordinarity holds.

Wan showed that Adolphson and Sperber's conjecture is in general false \cite{WAN93}.  Useing maximizing functions from linear programming, he obtained several decomposition theorems that, in effect, decompose the property of generic ordinarity \cite{WAN92}, \cite{WAN93},\cite{WAN04}.  

In this paper we show that the star and parallel hyperplane decomposition appearing in \cite{WAN93} and the collapsing decomposition that appears in \cite{WAN04} are each instances of a more general decomposition type referred to in linear programming as a regular decomposition.  This application of regular decompositions was first conjectured, by Wan in the case of toric hypersurfaces \cite{GOTT08}.  One must note the facial decomposition and boundary decomposition in \cite{WAN93} are not special cases of the regular decomposition.  In fact, the defining property of a regular decomposition is that it mimics a facial decomposition.  Both the facial and boundary decompositions are heavily used in showing that regular decomposition decomposes the property of ordinarity.

Throughout this work several examples of ordinarity, and several decompositions are provided.  We conclude with a demonstration of the regular decomposition in the case of Deligne polytopes.

\label{chapter_intro}

\subsection{Definition of $L$-function}
Let $p$ be a prime and $q=p^a$ for some positive integer $a$.  Let $\ff_q$ be the finite field of $q$ elements.  For each positive integer $k$, let $\ff_{q^k}$ be the finite extension of $\ff_q$ of degree $k$.  Let $\zeta_p$ be a fixed primitive $p$-th root of unity in the complex numbers.  For any Laurent polynomial $f(x_1,\ldots,x_n) \in \ff_q[x_1^{\pm1}, \ldots, x_n^{\pm1}]$, we form the exponential sum
\begin{equation}
\label{equationS}
S_k^*(f) = \sum_{x_i \in \ff_{q^k}^*} \zeta_p^{\Tr_{k}f(x_1,\ldots,x_n)},
\end{equation}
where $\ff_{q^k}^*$ denotes the set of non-zero elements in $\ff_{q^k}$ and $\Tr_k$ denotes the trace map from $\ff_{q^k}$ to the prime field $\ff_p$.  

By a theorem of Dwork-Bombieri-Grothendieck, the following generating $L$-function is a rational function \cite{DWORK62}, \cite{GROTHENDIECK64}:
\begin{equation} \label{Lfunction} L^*(f,T) = \exp \left( \sum_{k=1}^{\infty} S_k^*(f) \frac{T^k}{k} \right).\end{equation}

We may write $f$ as: 
$$f= \sum_{j=1}^J a_j x^{V_j}, a_j \neq 0,$$
where each $V_j = (v_{1j},\ldots, v_{nj})$ is a lattice point in $\zz^n$ and the power $x^{V_j}$ is the product $x_1^{v_{1j}} \cdots x_n^{v_{nj}}$.  Let $\Delta(f)$ be the convex closure in $\rr^n$ generated by the origin and the lattice points $V_j$ $(1 \leq j \leq J)$.  This is called the Newton polyhedron of $f$.  Without loss of generality we may always assume that $\Delta(f)$ is $n$-dimensional.  

For $\delta$ a subset of $\{V_1 ,\ldots, V_J\}$, we define the restriction of $f$ to $\delta$ to be the Laurent polynomial
$$f_{\delta} = \sum_{V_j \in \delta} a_j x^{V_j}.$$
For our purposes, we will generally take $\delta$ to be a sub-polytope or a face of $\Delta$.
This polytope structure suggests the following definition:
\begin{definition}
 The Laurent polynomial $f$ is called non-degenerate or $\Delta$-regular if for each closed face $\delta$ of $\Delta(f)$ of arbitrary dimension which does not contain the origin, the $n$ partial derivatives
$$\{ \frac{\partial f_{\delta}}{ \partial x_1 }, \ldots , \frac{\partial f_{\delta}}{ \partial x_n } \}$$
have no common zeros with $x_1 \cdots x_n \neq 0$ over the algebraic closure of $\ff_q$.
\end{definition}


In \cite{AS89}, Adolphson and Sperber proved that when $f$ is non-degenerate 
\begin{equation}
\label{eqnL}
L^*(f,T)^{(-1)^{n-1}} = \sum_{i=0}^{n!Vol(\Delta(f))} A_i(f) T^i, A_i(f) \in \zz[\zeta_p].
\end{equation}
In other words, $L^*(f,T)^{(-1)^{n-1}}$ is a polynomial of degree $n!\Vol(f)$.

\section{Newton Polygons and Hodge Polygons}
\subsection{Newton Polygons}
Let $g(x) = 1+ \sum_{i=1}^n a_i x^i \in 1 + x \zz[\zeta_p]$.  Consider the following sequence of points in the real coordinate plane:
$$(0,0), (1 ,\ord_q a_1), (2 ,\ord_q a_2),\ldots, (i ,\ord_q a_i),\ldots,(n ,\ord_q a_n),$$
where $\ord_q$ denotes the standard $q$-adic valuation on $\qq_p$, the field of $p$-adic numbers.  We normalize the valuation so that $\ord_q q = 1$.
If $a_i=0$ we omit that point.  Equivalently, we may think of it as lying ``infinitely" far above the horizontal axis since all finite powers of $p$ can be divided into $0$ without remainder.  The $q$-adic Newton polygon of $g(x)$ is defined to be the lower convex hull of this set of points, i.e. the highest convex polygonal line joining $(0,0)$ with $(n, \ord_q a_n)$ which passes on or below all of the points $(i, \ord_q a_i)$.  A more complete introduction to Newton polygons  and valuations appear in \cite{KOBLITZ1}.  The $q$-adic Newton polygon of (\ref{eqnL}) is denoted $NP(f)$.

Often it is convenient to think of $NP(f)$ as the real valued function on the interval $[0, n!\Vol(f)]$ whose graph is the Newton polygon.  Note that $L^*(f,T)^{(-1)^{n-1}}$ must be a polynomial and not just a rational function in order to define $NP(f)$.  Therefore restricting to the case where $f$ is non-degenerate guarantees that $NP(f)$ is well defined.  Since Newton polygons vary greatly as $f$ and $p$ vary, determining the Newton polygon is, in general, very complicated.  However, using Dwork theory we can obtain good estimates. 

For a fixed finite integral polytope $\Delta \subset \rr^n$, let $\mathcal{N}_p(\Delta)$ be the parameter space of $f$ over $\overline{\ff_p}$ such that $\Delta(f)=\Delta$.  This is a smooth affine variety defined over $\ff_p$.  Let $\mathcal{M}_p(\Delta)$ be the set of non-degenerate $f$ over $\overline{\ff_p}$ with $\Delta(f)=f$.  This is the compliment of a certain discriminant locus.  Thus, $\mathcal{M}_p(\Delta)$ is a Zariski open smooth affine subset of $\mathcal{N}_p(\Delta)$.  For $p$ sufficiently large, say $p > n!\Vol(\Delta)$, $\mathcal{M}_p(\Delta)$ is non-empty. 

For $f \in \mathcal{M}_p(\Delta)$, $NP(f)$ may vary greatly.   However, from the Grothendieck specialization theorem \cite{WAN00DWORK1} one may deduce that the lowest Newton polytope exists and is attained for all $f$ in some Zariski open dense subset of $\mathcal{M}_p(\Delta)$.  Hence we define the generic Newton polygon:
$$GNP(\Delta,p): = \inf_{f \in \mathcal{M}_p(\Delta)} NP(f).$$

Some work has been done to compute $GNP(\Delta,p)$ in the in the one-dimensional case \cite{ZHU03} and in certain two variable cases in \cite{NIU12}.
Though Newton polygons for specific $f$ may be out of our reach, we may be able to compute the lower bound of $NP(f)$ given by $GNP(\Delta,p)$.

\subsection{Definition of Hodge Polygon}

Newton polygons lie above a certain topological or combinatorial lower bound called the Hodge polygon.  This is given by Adolphson and Sperber in terms of rational points in $\Delta$.  This is notated $HP(\Delta)$.  Our construction of $HP(\Delta)$ will be strictly combinatoric.

For a given $f$, let $\Delta$ denote the $n$-dimensional integral polyhedron $\Delta(f)$ in $\rr^n$ containing the origin.  Let $C(\Delta)$ be the cone generated by $\Delta$ and the origin.  For a vector $u$ in $\rr^n$, $w(u)$ is defined to be the smallest positive real number $c$ such that $u \in c \Delta$.  If no such $c$ exists, that is, $u \notin C(\Delta)$, we define $w(u) = \infty$.

For $ u \in C(\Delta)$, the ray passing from the origin through $u$ intersects $\Delta$ in a face $\delta$ of co-dimension $1$.  This face is in general not unique unless the intersection point is in the interior of $\delta$.  Let $\sum_{i=1}^n e_i X_i=1$ be the equation of the hyperplane containing $\delta$ in $\rr^n$.  The coefficients $e_i$ are uniquely determined rational numbers.  One can show using linear programming the weight function is given by the formula
$$w(u) = \sum_{i=1}^{n} e_i u_i,$$
where $(u_1,\ldots,u_n)=u$ denotes the coordinates of $u$.

Let $D(\delta)$ be the least common denominator of the rational numbers $e_i$ for $1 \leq i \leq n$.  It follows that for a lattice point $u \in C(\delta)$, we have
$$w(u) \in \frac{1}{D(\delta)} \zz_{\geq 0}.$$
It can be shown that there are $u \in C(\delta)$ where $w(u)$ has denominator exactly $D(\delta)$.  In this sense $D(\delta)$ is optimal. Let $D(\Delta)$ be the least common denominator of all the $\delta$:
$$D(\Delta) = \lcm_{\delta} D(\delta),$$
where $\delta$ runs over all the co-dimension $1$ faces of $\Delta$ which do not contain the origin.  Thus we have 
$$w(\zz^n) \subseteq \frac{1}{D(\Delta)} \zz_{\geq 0} \cup \{ + \infty\}.$$

Let $D = D(\Delta)$.  For an integer $k$, let
$$W_{\Delta}(k) = card \{u \in \zz^n | w(u) = \frac{k}{D}\}.$$
This is the number of lattice points in $\zz^n$ with weight $k/D$.  Let 
$$H_{\Delta}(k) = \sum_{i=0}^{n} (-1)^i {n \choose i} W_{\Delta}(k-iD).$$ 
This number is the number of lattice points of weight $k/D$ in a certain fundamental domain corresponding to a basis of the $p$-adic cohomology space used to compute the $L$-function.  Thus $H_{\Delta}(k)$ is a non-negative integer for each $k\geq 0$.  Furthermore,
$$H_{\Delta}(k)=0, \mathrm{for~} k > nD$$
and 
$$\sum_{k=0}^{nD} H_{\Delta}(k)=n!\Vol(\Delta).$$
\begin{definition}
\label{defHP}
The Hodge polygon $HP(\Delta)$ of $\Delta$ is defined to be the lower convex polygon in $\rr^2$ with vertices
$$\left( \sum_{k=0}^{m} H_{\Delta}(k), \frac{1}{D} \sum_{k=0}^m k H_{\Delta}(k) \right).$$
That is, the polygon $HP(\Delta)$ is the polygon starting from the origin with a side of  slope $k/D$ with horizontal length $H_{\Delta}(k)$ for each integer $0 \leq k \leq nD$.
\end{definition}

As shown in \cite{AS96}, the numbers $H_{\Delta}(k)$ are the Hodge numbers in the toric hypersurface case, thus the term ``Hodge polygon."  

\subsection{Bounds on the Newton Polygon}
In \cite{AS89} Adolphson and Sperber proved that for $f \in \mathcal{M}_p(\Delta)$, 
$$NP(f) \geq HP(\Delta).$$
That is, the graph of $NP(f)$ lies above the graph of $HP(\Delta)$ at every point.  This is a Katz type conjecture.
From this we deduce:
\begin{prop}
For every prime $p$ and every $f \in \mathcal{M}_p(\Delta)$, we have the inequalities
$$NP(f) \geq GNP(\Delta,p) \geq HP(\Delta).$$
\end{prop}

\begin{definition}
If $NP(f) = HP(\Delta)$ then we say that $f$ is ordinary.
If $GNP(\Delta,p)= HP(\Delta)$, we say that the family $\mathcal{M}_p(\Delta)$ is generically ordinary.
\end{definition}

Adolphson and Sperber noticed the utility of $D(\Delta)$ in many Newton polygon computations. In \cite{AS87} they conjectured the following:
\begin{conjecture}[Adophson-Sperber]
\label{conjAS}
 If $ p\equiv 1 \modulo{D(\Delta)}$, then the $GNP(\Delta,p)$ coincides with  $HP(\Delta)$.
\end{conjecture}
This is a generalization of a conjecture of Dwork \cite[pg.40]{DWORK73} and Mazur \cite[pg.661]{MAZUR72}. Wan showed in \cite{WAN93} that Conjecture \ref{conjAS} is false for all $n\geq 5$.  However, he was able to weaken the conjecture proving:
\begin{theorem}[Wan]
\label{theoremD}
 There is an effectively computable integer $D^*(\Delta)$ such that if $p \equiv 1 \modulo{D^*(\Delta)}$, then $GNP(\Delta,p) = HP(\Delta)$.
\end{theorem}
Conjecture \ref{conjAS} is true in many important cases.  In general $D^*(\Delta)$ is difficult to compute.  

\subsection{Diagonal Laurent Polynomials}
A Laurent polynomial is called diagonal if $f$ has exactly $n$ non-constant terms and $\Delta(f)$ is $n$-dimensional.  The $L$-function can be computed explicitly using Gauss sums and the Stickelberger theorem.  They may also be used to show the following:
\begin{lemma}
For a diagonal Laurent polynomial
$$f = \sum_{j=1}^n a_j x^{V_j}, a_j \in \ff_q^*,$$
the Newton polygon depends only on the fact that the coefficients are nonzero. 
\label{lemmaDiag}
\end{lemma}
Thus in the diagonal case, for simplicity we may assume $f$ is of the form
\begin{equation}
f = \sum_{j=1}^n x^{V_j}, a_j \in \ff_q^*.
\label{eqDiagonal}
\end{equation}
Diagonal Laurent polynomials will provide important building blocks for many computations used in decomposition theory.

Let $M$ be the matrix of exponents of a diagonal Laurent polynomial as in (\ref{eqDiagonal}):
$$M= ( V_1,\ldots,V_n),$$
where each $V_j$ is written as a column vector.  One can check that $f$ is non-degenerate if and only if $p$ is relatively prime to $\det M$.  Consider the set of solutions to the following linear system:
\begin{equation}
M\left( \begin{matrix} r_1 \\ \vdots \\ r_n\end{matrix}\right) \equiv 0 \modulo{1}, r_i \mathrm{~rational~}, 0 \leq r_i <1.
\label{eqM}
\end{equation}
The map $(r_1,\ldots, r_n) \rightarrow r_1V_1+\ldots+ r_n V_n$ establishes a correspondence between the solutions to (\ref{eqM}) and the lattice points of the fundamental domain
$$\rr V_1+\ldots+\rr V_n \modulo{\zz V_1+\ldots \zz V_n}.$$
Let $S(\Delta)$ denote the set of solutions $r$ of (\ref{eqM}) which may be identified with the lattice points in the fundamental domain.  This has a natural abelian group structure under addition modulo $1$.  The order of $S(\Delta)$ is precisely given by 
$$\det M = n!\Vol(\Delta).$$
Let $S_p(\Delta)$ denote the prime to $p$ part of $S(\Delta)$.  It is an abelian subgroup of order equal to the prime to $p$ factor of $\det M$.  

\begin{theorem}
Suppose $f$ is a non-degenerate diagonal Laurent polynomial.  Also suppose that $\Delta=\Delta(f)$ and $p \nmid n!Vol(\Delta(f))$. Then $f$ is ordinary at $p$ if and only if the norm function $|r|$ on $S_p(\Delta)$ is stable under the $p$-action. That is, for each $r \in S_p(\Delta)$, we have
$$|r| = |\{pr\}|.$$
Where $|r| = r_1+\ldots+r_n$.  Equivalently, 
$$w(u) = w(\{pu\}).$$
In other words, the weight function $w(u)$ on the lattice points of $S_p(\Delta)$ is stable under the $p$-action.
\label{theoremDiag}
\end{theorem}
Thus we have established some conditions to detect ordinarity of Laurent polynomials.  The proofs of Lemma \ref{lemmaDiag} and Theorem \ref{theoremDiag} appear in \cite{WAN04}.

\subsection{An Example}

\begin{figure}[ht]
\begin{center}
\begin{tikzpicture}
\draw[to-to,very thin] (-0.5,0) -- (4.5,0); 
\draw[to-to,very thin] (0,-0.5) -- (0,3.5); 
\draw[step= 1cm,gray,very thin] (-0.25,-0.25) grid (4.25,3.25);
\draw[very thick] (0,0) -- (1,0) -- (4,3);
\draw (0,0) node {$\bullet$};
\draw (1,0) node {$\bullet$};
\draw (4,3) node {$\bullet$};
\end{tikzpicture}
\end{center}
\caption{The Hodge Polygon for $\frac{1}{x_1}+x_1x_2^2+x_1 x_3^2$}
\label{figure_ch1example}
\end{figure}
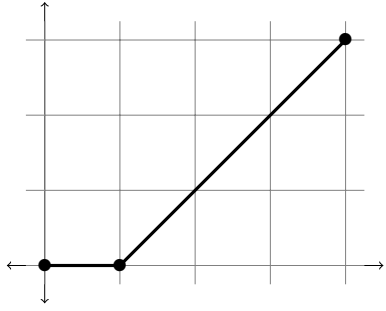

Let $p>2$ be prime and consider the polynomial 
$$f(x_1,x_2,x_3) = \frac{1}{x_1}+x_1x_2^2+x_1 x_3^2 \in \ff_p[x_1,x_2,x_3].$$  
The polytope $\Delta(f)=\Delta$ is spanned by the origin and the vertices $\smallpmatrix{-1 \\ 0 \\ 0}, \smallpmatrix{1 \\2 \\0}$, and $\smallpmatrix{1\\0\\2}$.  It follows that $D(\Delta) = 1$.  The Hodge polygon is computed by first finding the values for $W_{\Delta}(k)$ and $H_{\Delta}(k)$.  These values are summarized in Table \ref{tableWandH}.  Using (\ref{eqnL}) we need only to check exponents up to $nD=3$.  
\begin{table}[ht]
\caption{$W_{\Delta}(k)$ and $H_{\Delta}(k)$}
\centering
\begin{tabular}{c | c c c c }
\hline\hline
k & 0& 1& 2& 3\\ \hline
$W_{\Delta}(k)$ & 1 & 6 & 15 & 28 \\
$H_{\Delta}(k)$ & 1 & 3 & 0 & 0\\
[1 ex] \hline
\end{tabular}
\label{tableWandH}
\end{table}
In general:
$$W_{\Delta}(k) = (k+1)(k+2)/2.$$
We can generate the vertices of $HP(\Delta)$ using Definition \ref{defHP}.  This information is summarized Table \ref{tableHPvertices}.  Figure \ref{figure_ch1example} is the graph of $HP(\Delta)$.

\begin{table}[h]
\caption{Hodge Polygon Vertices}
\centering
\begin{tabular}{c | c c c c}
\hline\hline
m & 0& 1& 2& 3\\ \hline
$\sum_{k=0}^m H_{\Delta}(k)$ & 1 &1 & 4& 4\\
$\frac{1}{D}\sum_{k=0}^m k H_{\Delta}(k)$ & 0 & 0 & 3& 3 \\
[1 ex] \hline
\end{tabular}
\label{tableHPvertices}
\end{table}

For $p=q=3$, we may compute the polynomial $L(f,T)^{(-1)^{n-1}}=L(f,T)$ directly.  In general this is computationally very difficult.  However, for small primes and an efficient computer algebra package, such calculations are possible.  Via direct computation we determine the $S_{\Delta}(k)$ values summarized in Table \ref{tableS}.  Since $f$ is non-degenerate for $p>2$, we know that $L(f,T)$ is a polynomial of degree $n!\Vol(\Delta) = |\det\smallpmatrix{-1& 1 & 1 \\ 0 & 2 & 0 \\ 0& 0 & 2}|=4$.
From this we determine that 
$$L(f,T) = -27T^4+18T^2+8T+1.$$  The vertices of the Newton polygon we generate are $(0,0) ,(1,0), (2,2)$ and $(4,3)$.  The lower convex hull of these points coincide with $HP(\Delta)$.  Therefore $NP(f)=HP(\Delta(f))$ and $f$ is ordinary.

Note that at $p=2$ so the Newton polygon of $f$ is not defined.  One can also show that $f$ is ordinary for all primes other than 2.  Hence $D^*(\Delta)=1$.  For $p=q=2$ through direct computation we find that $L$-function is $1-T$.  

A more robust generalization of this example, which include conditions when $\MM_p(\Delta)$ is not generically ordinary is given in Section \ref{chapter_deligne}.  

\begin{table}[h]
\caption{Exponential Sum Values}
\centering
\begin{tabular}{c | c c c c}
\hline\hline
$k$ & 1 & 2 & 3 & 4\\ \hline
$S_{\Delta}(k)$ & 8 & -14 & 80/3 & -61\\ [1 ex] \hline
\end{tabular}
\label{tableS}
\end{table}
 
The calculation of the $L$-function is largely dependent on the existence of the terms, rather than the individual values of the coefficient.  Also recall the Hodge polygon was based solely on the Newton polyhedra $\Delta$.  This suggests a strong connection between the geometry of $\Delta$ and the shape of the Newton polygon.  Wan showed that certain decompositions of $\Delta$ induce a decomposition of the property of ordinarity.  We devote the remainder of this work to highlight such decompositions and establish a new decomposition, Theorem \ref{mainth}, the Regular Decomposition theorem.

\section{Dwork Theory}
\label{chapter_dwork}
To examine ordinarity we must first examine $L(f,T)$ in more depth.  Dwork's trace formula will allow us to express the $L$-function as the Fredholm determinant of a certain infinite Frobenius matrix.  From this we will descend to a related infinite matrix which captures the behavior of ordinarity on $L(f,T)$.  Much of this development is taken from \cite[\S 4]{WAN04}.  In all of our work we will assume $f$ is non-degenerate.

\subsection{Dwork's Trace Formula}
As before, let $q=p^a$ where $p$ is prime and $a$ is some positive integer.  Let $K$ denote the unramified extension of $\qq_p$ in $\Omega$ of degree $a$.  Let $\Omega_1 = \qq_p(\zeta_p)$ where, as before, $\zeta_p$ denotes a primitive $p$-th root of unity.  Thus $\Omega_1$ is the totally ramified extension of $\qq_p$ of degree $p-1$.  Let $\Omega_a$ be the compositum of $ \Omega_1$ and $K$.  The field $\Omega_a$ is an unramified extension of $\Omega_1$ of degree $a$.  The residue fields of $\Omega_a$ and $K$ are both $\ff_q$, and the residue fields of $\Omega_1$ and $\qq_p$ are both $\ff_p$.  Let $\pi$ be a fixed element in $\Omega_1$ satisfying
$$\sum_{m=0}^{\infty} \frac{\pi^{p^m}}{p^m} = 0, \ord \pi = \frac{1}{p-1}.$$
Hence $\pi$ is a uniformizer of $\Omega_1 = \qq_p(\zeta_p)$ and we have
$$\Omega_1 = \qq_p(\pi).$$
The Frobenius automorphism $x \mapsto x^p$ of $\Gal(\ff_q/ \ff_p)$ lifts to a generator $\tau$ of $\Gal(K/\qq_p)$ which is extended to $\Omega_a$ by requiring that $\tau(\pi) = \pi$.  If $\zeta$ is a $(q-1)$-st root of unity in $\Omega_a$, then $\tau(\zeta)= \zeta^p$.

Let $E(t)$ be the Artin-Hasse exponential series:
$$
\begin{array}{rcl}
E(t) & = & \exp(\sum_{m=0}^{\infty} \frac{t^{p^m}}{p^m})\\
      & = & \prod_{k \geq 1, (k,p)=1} (1-t^k)^{\mu(k)/k},
\end{array}$$
where $\mu(k)$ is the M\"obius function.  The last product expansion shows that the power series $E(t)$ has $p$-adic integral coefficients.  Thus we can write
$$E(t) = \sum_{m=0}^{\infty} \lambda_m t^m, \lambda_m \in \zz_p.$$
For $0 \leq m \leq p-1$, more precise information is given by
$$\lambda_m = \frac{1}{m!}, \ord \lambda_m = 0, 0 \leq m \leq p-1.$$
The shifted series
$$\theta(t) = E(\pi t) = \sum_{m=0}^{\infty} \lambda_m \pi^m t^m$$
is a splitting function in Dwork's terminology.  The value $\theta(1)$ is a primitive $p$-th root of unity, which will be identified with the $p$-th root of unity $\zeta_p$ used in our definition of the exponential sum in (\ref{equationS}).

For a Laurent polynomial $f(x_1,\ldots,x_n) \in \ff_q[x_1, x_1^{-1},\ldots, x_n, x_n^{-1}]$, we write
$$f= \sum_{j=1}^J \overline{a}_j x^{V_j}, V_j \in \zz^n, \overline{a}_j \in \ff_q^*.$$
Let $a_j$ be the Teichm\"uller lifting of $\overline{a}_j$ in $\Omega$.  Thus, we have $a_j^q=a_j$.  Set 
$$F(f,x) = \prod_{j=1}^J \theta(a_j x^{V_j})$$
$$F_a(f,x) = \prod_{i=0}^{a-1} F^{\tau^i}(f,x^{p^i}).$$
Note that by the definition of $\theta$, $F(f,x)$ and $F_a(f,x)$ are well defined as formal Laurent series in $x_1,\ldots,x_n$ with coefficients in $\Omega_a$.

To describe the growth conditions satisfied by $F$, write
$$F(f,x) = \sum_{r \in \zz^n} F_r(f) x^r.$$
Then from the definition of $F$ and $F_a$ one checks that
\begin{equation}
\label{eq45}
F_r(f) = \sum_u \left( \prod_{j=1}^J \lambda_{u_j} a_j^{u_j}  \right) \pi^{u_1+\ldots+u_J},
\end{equation}
where the outer sum is over all solutions of the linear system
\begin{equation}
\label{eq46}
\sum_{j=1}^J u_j V_j =r , u_j \geq 0, u_j \mathrm{~integral}.
\end{equation}
Thus $F_r(f)=0$ if (\ref{eq46}) has no solutions.  Otherwise (\ref{eq45}) implies that
\begin{equation}
\ord F_r(f) \geq \frac{1}{p-1} \inf_u \{\sum_{j=1}^J u_j\},
\label{eq46b}
\end{equation}
where the infimum is taken over all solutions of (\ref{eq46}).

For $r\ \in \rr^n$, recall that the weight $w(r)$ is given by
$$w(r) = \inf_u \{ \sum_{j=1}^J u_j \mid \sum_{j=1}^J u_j v_j = r, u_j \geq 0, u_j \in \rr \},$$
where $w(r)$ is taken to be $\infty$ if $r$ is not in the cone generated by $\Delta$ and the origin.  Thus
\begin{equation}
\label{eq47}
\ord F_r(f) \geq \frac{w(r)}{p-1},
\end{equation}
with the convention that $F_r(f) = 0$ if $ w(r) = +\infty$.

Recall we defined $C(\Delta)$ to be the closed cone generated by the origin and $\Delta$.  Let $L(\Delta)$ be the set of lattice points in $C(\Delta)$.  That is,
$$L(\Delta) = \zz^n \cap C(\Delta).$$
For real numbers $b$ and $c$ with $0 \leq b \leq p/(p-1)$, define the following two spaces of $p$-adic functions:
$$\frL(b,c) =\{ \sum_{r \in L(\Delta)} C_r x^r \mid C_r \in \Omega_a, \ord_p C_r \geq bw(r)+c \}$$
$$\frL(b) = \bigcup_{c \in \rr} \frL(b,c).$$
One checks from (\ref{eq47}) that
$$F(f,x) \in \frL(\frac{1}{p-1},0), ~~F_a(f,x) \in \frL(\frac{p}{q(p-1)},0).$$
Define an operator $\psi$ on formal Laurent series by
$$\psi( \sum_{r \in L(\Delta)} C_r x^r) = \sum_{r \in L(\Delta)} C_{pr} x^r.$$

Therefore 
$$\psi(\frL(b,c)) \subset \frL(pb,c).$$
It follows that the composite operator $\phi_a = \psi^a \circ F_a(f,x)$ is an $\Omega_a$-linear endomorphism of the space $\frL(b)$, where $F_a(f,x)$ denotes the multiplication map by the power series $F_a(f,x)$.  Similarly, the operator $\psi_a = \tau^{-1} \circ \psi \circ F(f,x)$ is an $\Omega_a$-semilinear ($\tau^{-1}$-linear) endomorphism of the space $\frL(b)$.  The operators $\phi_a^m$ and $\phi_1^m$ have well defined traces and Fredholm determinants.  The Dwork trace formula asserts that for each positive integer $k$,
$$S_k^*(f) = (q^k-1)^n\Tr(\phi_a^k).$$
Applying this to $L^*(f,t)^{(-1)^{n-1}}$ we have:
\begin{theorem}
$$L^*(f,T)^{(-1)^{n-1}} = \prod_{i=0}^n \det(I-Tq^i \phi_a)^{(-1)^i { n \choose i}}.$$
\end{theorem}

Hence, understanding of the $L$-function is reduced to understanding the single determinant $\det(I-T\phi_a)$.  For a more tangible representation we shall describe the operator $\phi_a$ in terms of an infinite nuclear matrix.  First, observe that 
$$\begin{array}{rcl}
\phi_1^a & = & \phi_1^{a-2} \tau^{-1} \circ \psi \circ F(f,x) \circ \tau^{-1} \circ \psi \circ F(f,x) \\
               & = & \phi_1^{a-2} \tau^{-1} \circ \tau^{-1} \circ F^{\tau}(f,x) \circ \psi \circ F(f,x) \\
               & = & \phi_1^{a-2} (\tau^{-2}) \circ \psi^2 \circ F^{\tau}(f,x^p) \circ F(f,x) \\
               & = & \ldots = \psi^a \circ F_a(f,x) = \phi_a.
\end{array}$$
We now describe the matrix form of the operators $\phi_1$ and $\phi_a$ with respect to some orthonormal basis. Let $D=D(\Delta)$ as defined in Section \ref{chapter_intro}.  Fix a choice $\pi^{1/D}$ of $D$-th root of $\pi$ in $\Omega$.  The monomials $\pi^{w(r)} x^r$ (where $r \in L(\Delta))$ form an orthonormal basis of the $p$-adic Banach space 
$$\mathcal{B}: = \{ \sum_{r \in L(\Delta)} C_r \pi^{w(r)} x^r \mid C_r \in \Omega_a(\pi^{1/D}), C_r \rightarrow 0 \mathrm{~as~}|r| \rightarrow \infty\}.$$
Furthermore, if $b> 1/(p-1)$, then $\frL(b) \subseteq \mathcal{B}$.  The operator $\phi_a$ (resp. $\phi_1$) is an $\Omega_a$-linear (resp. $\Omega_a$-semilinear) nuclear endomorphism of the space $\mathcal{B}$.  Let $\Gamma$ be the orthonormal basis $\{\pi^{w(r)}x^r\}_{r\in L(\Delta)}$ of $\mathcal{B}$ written as a column vector.  One checks that the operator $\phi_1$ is given by
$$\phi_1 \Gamma = A_1(f)^{\tau^{-1}}\Gamma,$$
where $A_1(f)$ is the infinite matrix whose rows are indexed by $r$ and columns are indexed by $s$.  That is,
\begin{equation}
\label{eq49}
A_1(f) = (a_{r,s}(f)) = (F_{ps-r}(f)\pi^{w(r)-w(s)}).
\end{equation}
Note that the $\pi^{w(r)-w(s)}$ factor of each term is the contribution to the $p$-adic valuation that guarantees $NP(f) \geq HP(\Delta)$.

Since $\phi_a = \phi_1^a$ and $\phi_1$ is $\tau^{-1}$-linear, the operator $\phi_a$ is given by
$$\begin{array}{rcl}
\phi_a \Gamma & = & \phi_1^a \Gamma \\
			 & = & \phi_1^{a-1} A_1^{\tau^{-1}} \Gamma \\
			 & = & \phi_1^{a-2} A_1^{\tau^{-2}}A_1^{\tau^{-1}} \Gamma \\
			 & = & A_1^{\tau^{-a}} \ldots A_1^{\tau^{-2}}A_1^{\tau^{-1}} \Gamma \\
			 & = & A_1 A_1^{\tau} \ldots A_1^{\tau^{a-2}}A_1^{\tau^{a-1}} \Gamma. 
\end{array}$$
Let
$$A_a(f) = A_1 A_1^{\tau^1} \ldots A_1^{\tau^{a-1}}.$$
Then the matrix of $\phi_a$ under the basis $\Gamma$ is $A_a(f)$.  We call $A_1(f) = (a_{r,s}(f))$ the infinite semilinear Frobenius matrix and $A_a(f)$ the infinite linear Frobenius matrix.  Dwork's trace formula can now be rewritten in terms of the matrix $A_a(f)$ as follows:
\begin{equation}
\label{eq51}
L^*(f,T)^{(-1)^{n-1}} = \prod_{i=0}^n \det(I-Tq^iA_a(f))^{(-1)^i { n \choose i}}.
\end{equation}
Hence we are now reduced to understanding the single determinant $\det(I-TA_a(f))$.

\subsubsection{Newton Polygons of Fredholm Determinants}
To get a lower bound for the Newton polygon of $\det(I-TA_a(f))$, we need to estimate the entries of the infinite matrices $A_1(f)$ and $A_a(f)$.  By (\ref{eq47}) and (\ref{eq49}), we obtain the estimate
\begin{equation}
\label{estimate52}
\ord a_{r,s}(f) \geq \frac{w(ps-r)+w(r)-w(s)}{p-1}\geq w(s).
\end{equation}
Recall that for a positive integer $k$, $W_{\Delta}(k)$ is defined to be the number of lattice points in $L(\Delta)$ with weight exactly $k/D$:
\begin{equation}
\label{eqW}
W_{\Delta}(k) = \#\{r \in L(\Delta) \mid w(r) = \frac{k}{D}\}.
\end{equation}
Let $\xi \in \Omega$ such that 
$$\xi^{D} = \pi^{p-1}.$$
Therefore $\ord \xi = 1/D$.
By (\ref{estimate52}) the infinite matrix $A_1(f)$ has the block form
\begin{equation}
\label{block53}
A_1(f) = \left( \begin{array}{ccccc}
A_{00} & \xi^1 A_{01} & \ldots & \xi^i A_{0i} & \ldots \\
A_{10} & \xi^1 A_{11} & \ldots & \xi^iA_{1i} & \ldots \\
\vdots & \vdots & \ddots & \vdots &\\
A_{i0} & \xi^1 A_{i1} & \ldots & \xi^iA_{ii} & \ldots \\
\vdots & \vdots & \ddots & \vdots &
\end{array} \right)
\end{equation}
where the block $A_{ij}$ is a finite matrix of $W_{\Delta}(i)$ rows and $W_{\Delta}(j)$ columns whose entries are $p$-adic integers in $\Omega$.  The $\xi^i$ factors are the collection the $\pi^{w(r)-w(s)}$ terms in (\ref{eq49}).

\begin{definition}
Let $P(\Delta)$ be the polygon in $\rr^2$ with vertices $(0,0)$ and
$$\left( \sum_{k=0}^m W_{\Delta}(k), \frac{1}{D(\Delta)} \sum_{k=0}^m k W_{\Delta} (k) \right), ~m=0,1,2,\ldots.$$
\end{definition}
This is the chain level version of the Hodge polygon.  The block form in (\ref{block53}) and the standard determinant expansion of the Fredholm determinant show that we have:
\begin{prop}
The Newton polygon of $\det(I-TA_1(f))$ computed with respect to $p$ lies above the polygon $P(\Delta)$.
\end{prop}
Note that weighing the sum by $k$ in the second coordinate is a direct consequence of the $\xi^k$ factors of the block form matrix in (\ref{block53}).

Using the block form (\ref{block53}) and the exterior power construction of a semi-linear operator, one then gets the following lower bound of Adolphson and Sperber \cite{AS89} for the Newton polygon of $\det(I-TA_a(f))$.

\begin{prop}
The Newton polygon of $\det(I-TA_a(f))$ computed with respect to $q (=p^a)$ lies above the polygon $P(\Delta)$.
\end{prop}

\subsubsection{A Descent Theorem}
In general, the Newton polygon of $\det(I-T A_a(f))$ computed with respect to $q$ is different from the Newton polygon of $\det(I-TA_1(f))$ computed with respect to $p$, even though they have the same lower bound.  Since the matrix $A_a(f)$ is much more complicated than $A_1(f)$, especially for large $a$, we would like to replace $A_a(f)$ by the simpler matrix $A_1(f)$.  This is not possible in general.  However, if we are only interested in the question of whether the Newton polygon of $\det(I-TA_a(f))$ coincides with its lower bound, the following theorem shows that we can descend to the simpler $\det(I-TA_1(f))$.  


We are able to reduce ordinarity in terms of chain level polytopes.  If the Newton polygon of $\det(I-TA_1(f))=P(\Delta)$ we say that $f$ (or $\det(I-TA_1(f))$) is chain level ordinary.  The connection between chain-level ordinarity and generic ordinarity is established in \cite{WAN04}:

\begin{theorem}
\label{th4.7}
Let $\Delta(f) = \Delta$.  Assume that the $L$-function $L^*(f,T)^{(-1)^{n-1}}$ is a polynomial.  Then $NP(f) = HP(\Delta)$ if and only if the Newton polygon of $\det(I-TA_1(f))$ coincides with its lower bound $P(\Delta)$.  In this case the degree of the polynomial $L^*(f,T)^{(-1)^{n-1}}$ is exactly $n!V(f)$.
\end{theorem}

Hence our study of ordinarity is reduced to examining the ordinarity of $\det(I-TA_1(f))$.

\subsection{Boundary and Facial Decomposition Theorems}
To establish a regular decomposition theory, we must first examine the impact of decompositions of $\Delta$ on $\det(I-TA_1(f))$.  We present the first two decompositions here.  More work will is necessary before we introduce the others.

Let $B(\Delta)$ be the unique interior decomposition of the cone $C(\Delta)$ into a union of disjoint, relatively open cones.  Its elements are the interiors of the closed faces in $C(\Delta)$ that contain the origin.  The interior of the cone is the unique element in $B(\Delta)$ of dimension $n$.  The origin itself (if it is a vertex) is the unique element of $B(\Delta)$ of dimension $0$.  For $\Sigma \in B(\Delta)$, let $A_1(\Sigma,f)$ be the ``$\Sigma$" piece of $(a_{s,r}(f))$ in $A_1(f)$ i.e., $r$ and $s$ run through the cone $\Sigma$ rather than all of $C(\Delta)$.  In particular, for the full cone $\Sigma=C(\Delta)$, we have
$$A_1(C(\Delta),f) = A_1(f).$$
Let $A_1(\Sigma, f_{\Sigma})$ be the ``interior piece" of the Frobenius matrix $A_1(f_{\Sigma})$, where $f_{\Sigma}$ is the restriction of $f$ to the closure of $\Sigma$.

\begin{theorem}[Boundary Decomposition (Wan, \cite{WAN93})]
\label{th5.1}
The following factorization is true:
\begin{equation}
\label{eq5.1}
\det(I-TA_1(f)) = \prod_{\Sigma \in B(\Delta)} \det(I-TA_1(\Sigma, f_{\Sigma})).
\end{equation}
\end{theorem}
\begin{cor}
Let $\Sigma \in B(\Delta)$.  If the Newton polygon of $\det(I-TA_1(f))$ coincides with $P(\Delta)$, then the Newton polygon of $\det(I-TA_1(f_{\Sigma}))$ coincides with its lower bound $P(\overline{\Sigma})$.
\end{cor}

\begin{figure}[t]
\begin{center}

\begin{tikzpicture}
\draw[to-to,very thin] (-0.5,0) -- (4.5,0); 
\draw[to-to,very thin] (0,-0.5) -- (0,4.5); 
\draw[very thick] (0,0) -- (1,4) --  (4,1) -- (0,0);
\draw[step= 1cm,gray,very thin] (-0.25,-0.25) grid (4.25,4.25);
\end{tikzpicture}
~
\begin{tikzpicture}
\draw[to-to,very thin] (-0.5,0) -- (4.5,0); 
\draw[to-to,very thin] (0,-0.5) -- (0,4.5); 
\fill[black!10!white] (0,0) -- (1.0625,4.25) -- (4.25,4.25) -- (4.25,1.0625) -- (0,0);
\draw[step= 1cm, thin, gray] (-0.25,-0.25) grid (4.25,4.25);
\draw[-to,very thick, dashed] (0.018,0.072) -- (1.0625,4.25);
\draw[-to,very thick, dashed] (0.072,0.018) -- (4.25,1.0625);
\draw[dotted] (1,4) --(4,1);
\draw (0,0) node {$\circ$};
\end{tikzpicture}
\end{center}
\caption{Boundary Decomposition of $C(\Delta)$}
\label{figureBoundary}
\end{figure}
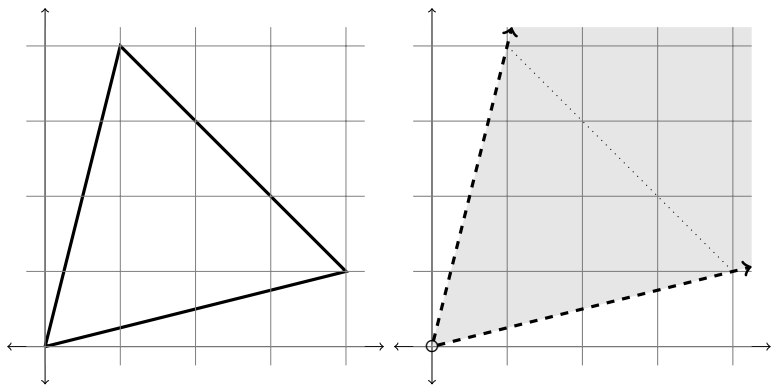

Figure \ref{figureBoundary} shows how a cone $C(\Delta)$ generated by the origin and the vertices $(4,1)$ and $(1,4)$ is decomposed into four sub-cones under the boundary decomposition:  the origin, the ray emanating from the origin and passing through $(4,1)$, the ray emanating from the origin and passing through $(1,4)$ and the open interior of $C(\Delta)$.

The next decomposition theorem, the facial decomposition theorem, was first obtained in \cite{WAN92}.   Let $\sigma_1,\ldots,\sigma_h$ be the $(n-1)$ dimension closed faces of $\Delta$, which do not contain the origin.  For each $1 \leq i \leq h$, let $f_{\sigma}$ be the restriction of $f$ to the closed polyhedron generated by $\sigma_i$ and the origin.  If $f$ is nondegenerate, then $f_{\sigma_i}$ is nondegenerate.  

\begin{theorem}[Facial Decomposition Theorem]
Let a Laurent polynomial $f$ be non-degenerate and let $\Delta(f)$ be $n$-dimensional.  Then $f$ is ordinary if and only if each $f_{\sigma_i}$ is ordinary.  Equivalently, $f$ is non-ordinary if and only if if some $f_{\sigma_i}$ is non-ordinary.
\label{theoremFacial}
\end{theorem}

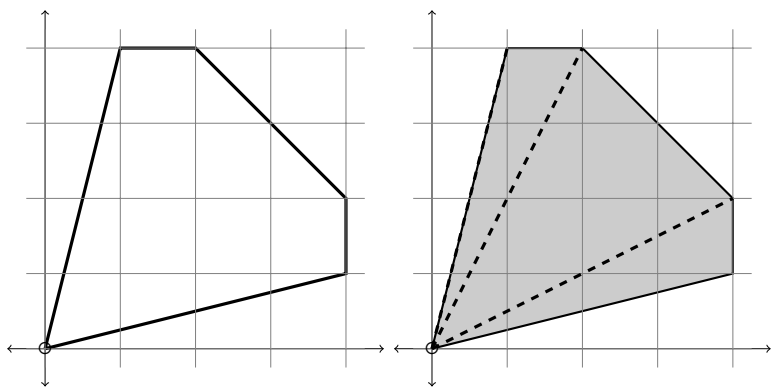
\begin{figure}[t]
\begin{center}
\begin{tikzpicture}
\draw[to-to,very thin] (-0.5,0) -- (4.5,0); 
\draw[to-to,very thin] (0,-0.5) -- (0,4.5); 
\draw[very thick] (0,0) -- (1,4) -- (2,4) -- (4,2) -- (4,1) -- (0,0);
\draw[step= 1cm,gray,very thin] (-0.25,-0.25) grid (4.25,4.25);
\draw (0,0) node {$\circ$};
\end{tikzpicture}
~
\begin{tikzpicture}
\draw[to-to,very thin] (-0.5,0) -- (4.5,0); 
\draw[to-to,very thin] (0,-0.5) -- (0,4.5); 
\fill[black!20!white] (0,0) -- (1,4) -- (2,4) -- (4,2) -- (4,1) -- (0,0);
\draw[thick] (0,0) -- (1,4) -- (2,4) -- (4,2) -- (4,1) -- (0,0);
\draw[step= 1cm,gray,very thin] (-0.25,-0.25) grid (4.25,4.25);
\draw[very thick, dashed] (0,0) -- (1,4);
\draw[very thick, dashed] (0,0) -- (2,4);
\draw[very thick, dashed] (0,0) -- (4,2);
\draw (0,0) node {$\circ$};
\end{tikzpicture}
\end{center}
\caption{Facial Decomposition}
\label{figureFacial}
\end{figure}

Figure \ref{figureFacial} shows the facial decomposition of a polytope with three faces not containing the origin.  

Theorem \ref{theoremFacial} allows us to assume that $\Delta$ is generated by a single face $\delta$ not containing the origin.  Using this we may assume that $\Delta$ contains a unique face not containing the origin.  We can now turn our focus toward decomposing $\delta$ which will in turn decompose the cone $C(\Delta)$.


\section{Polytope and Polygonal Constructions}
\label{chapter_poly}
\subsection{Decompositions}
\label{section_Decomposition}

\begin{figure}[t]

\begin{center}
\begin{tikzpicture}
\fill[black!20!white] (1,-0.45) -- (1,1) -- (2,3) -- (4,4) -- (5,2) -- (5,-0.45);
\draw (1,-0.5) -- (1,1) -- (2,3) -- (4,4) -- (5,2) -- (5,-0.5);
\draw (2,3) -- (2,-0.5);
\draw (4,4) -- (4,-0.5);
\draw[to-to] (-1,0) -- (6,0); 
\draw[to-to] (0,-1) -- (0,4.5); 
\draw (6.3,0) node{$\rr^n$};
\draw (0,4.8) node{$\rr$};
\draw[dashed] (4,4) -- (0,4);
\draw (4.2,0.2) node {$\omega$};
\draw (4,0) node {$\bullet$};
\draw (-0.5,4) node {$\psi(\omega)$};
\draw[dashed] (1.5,1) -- (0,1);
\draw (1.5,1) -- (1.5,-0.5);
\draw (1.75,0.25) node {$\omega'$};
\draw (1.5,0) node {$\bullet$};
\draw (-0.5,1) node {$\psi(\omega')$};
\draw (5.2,-0.3) node {$\delta$};
\end{tikzpicture}
\end{center}
\caption{Graph of $G_{\psi}$}
\label{figure_triangulation_example}
\end{figure}
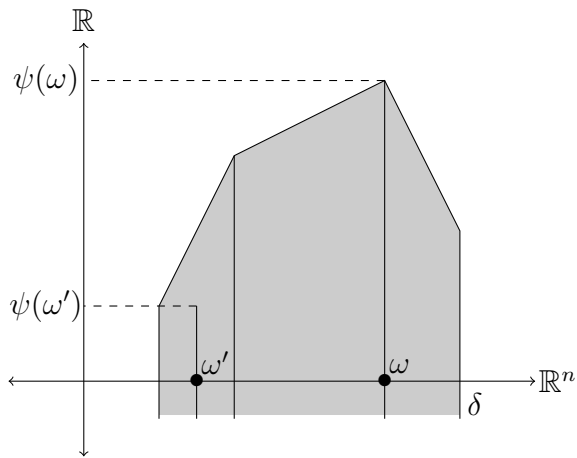

The facial decomposition allows us to assume that $\Delta$ contains a unique face $\delta$ that does not contain the origin.  In fact we will be decomposing $\delta$ and considering sub-cones of $C(\Delta)$ induced by this decomposition.

In general, decompositions are not required to have integral vertices, however for the purposes of examining $L$-functions it is more useful to restrict our discussion to polytopes with integral vertices and decompositions into integral polytopes.  We will assume all decompositions are integral.  A more general discussion of polytope decompositions can be found in \cite{GKZ94}.

Not all integral decompositions will induce a decomposition of ordinarity as in the facial decomposition.  To understand decomposition theorems more fully in the context of ordinarity on the chain level, we must first introduce some concepts from polytope decomposition theory and linear programming.  

\begin{definition}
A decomposition $T=\{\delta_1,\ldots,\delta_j\}$ of a convex polytope $\delta \subset \rr^n$ is a finite collection of $n$-dimensional polytopes (not necessarily convex) where
$$\delta = \delta_1 \cup \delta_2 \cup \ldots \cup \delta_h$$
and each $1\leq i \neq j \leq h$, $\delta_i \cap \delta_j$ is at worst a polytope of dimension $(n-1)$.  
\end{definition}
Note that the boundary decomposition is not a decomposition under this definition.  Neither is the facial decomposition.  We will always specify when we are using these two special decompositions.  
\subsection{Triangulations}
A triangulation of a convex polytope $\delta \subset \rr^n$ is a decomposition of $\delta$ into a finite number of simplices such that the intersection of any two of these simplices is a common face of them both (possibly empty).  Notationally we regard a triangulation as a collection of its simplices of maximal dimension.  All the lower-dimensional simplices are just faces of maximal ones.  Formally this is stated as follows:
\begin{definition}
 A triangulation $T=\{\delta_1,\ldots,\delta_h\}$ of $\delta$ is a decomposition where each $\delta_i$ is a simplex and $\delta_i \cap \delta_j$ is a common face for both $\delta_i$ and $\delta_j$.
\end{definition}
Triangulations are one of the most common and intuitive types of decompositions.

We are interested in decompositions whose vertices belong to a fixed finite set of lattice points.  Let $A$ be a finite subset of $\delta$ containing all the vertices of $\delta$ (therefore $\delta$ is the convex closure of $A$).  By a triangulation of $(\delta,A)$, we mean a triangulation of $\delta$ into simplices with vertices in $A$.  Note that we do not require every element of $A$ to appear as a vertex of a simplex.

\subsubsection{Construction of Decompositions}

Suppose $\delta$ is dimension $k$ and $A= \{V_1,\ldots,V_h\}$ is a subset of integral points which contain all vertices of $\delta$.
Take any function $\psi:A \rightarrow \rr$ and consider, in the space $\rr^{k+1} = \rr^k \times \rr$, the union of vertical half-lines
$$\{(\omega,y) : y \leq \psi(\omega), \omega \in A, y \in \rr\}.$$
Let $G_{\psi}$ be the convex hull of all these half lines (Figure \ref{figure_triangulation_example}).  This is an unbounded polyhedron projecting onto $\delta$.  The faces of $G_{\psi}$ which do not contain vertical half-lines (i.e. are bounded) form the bounded part of the boundary of $G_{\psi}$, which we call the upper boundary of $G_{\psi}$.  Clearly, the upper boundary projects bijectively onto $\delta$.  
If the function $\psi$ is chosen to be suitably generic, then all the bounded faces of $G_{\psi}$ are simplices and therefore their projections to $\delta$ define a unique triangulation of $(\delta,A)$.  However, for our work,we will only need decompositions, not triangulations.

Let $T$ be an arbitrary decomposition of $(\delta,A)$ where every element of $A$ is a vertex of some member of $T$.  Let $\psi:A \rightarrow \rr$ be any function.  Then there is a unique piecewise linear function $g_{\psi}:\delta \rightarrow \rr$ such that, for each $\omega \in A$ we have $g_{\psi}(\omega) = \psi(\omega)$.  The function $g_{\psi}$ is obtained by affinely interpolating $\psi$ inside each sub-polytope.  Each $\psi$ generates a unique $g_{\psi}$.  Each $g_{\psi}$ is a $T$-piecewise linear function.  For the purposes of generic ordinarity, we would like $\psi$ to be chosen sufficiently generic so that the domains of linearity are precisely the members of $T$.  This is not always possible and will lead to certain limitations in decomposition theory.  This construction will be useful later when we discuss maximizing functions.

\begin{figure}[t]
\begin{center}
\begin{tikzpicture}
\draw (0,0) -- (3,0) -- (1.5,2.6) -- (0,0);
\draw (0.75,0.433) -- (2.25,0.433) -- (1.5, 1.733) -- (0.75,0.433); 
\draw (0,0) -- (0.75,0.433);
\draw (3,0) -- (2.25,0.433);
\draw (1.5,2.6) -- (1.5, 1.733);
\end{tikzpicture}
\end{center}
\caption{A Regular Decomposition}
\label{figure_convex_triangulation}
\end{figure}
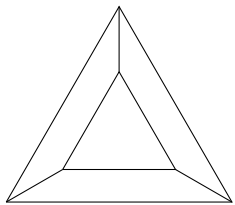

\subsubsection{Regular Decompositions}

As mentioned previously, triangulations are often the most useful in practice.  However other forms of decomposition are also useful.  In particular, Theorem \ref{mainth}, our main theorem, uses a regular decomposition that is not necessarily a triangulation. 

We are now able to define the regular decompositions referred to in Theorem \ref{mainth}.  A regular decomposition of $\delta$ is a decomposition $T$
into polytopes $\delta_1,\ldots,\delta_h$ such that  there is a piecewise linear function $\phi:\delta \mapsto \rr$ such that
        \begin{enumerate}
            \item $\phi$ is concave i.e. $\phi(t x+ (1-t) x') \geq t \phi(x) + (1-t) \phi(x')$, for all $x,x' \in \delta, 0\leq t \leq 1 $.
            \item The domains of linearity of $\phi$ are precisely the $\delta_i$ for $1 \leq i \leq m$.
        \end{enumerate}

\begin{figure}[t]
\begin{center}
\begin{tikzpicture}
\draw (0,0) -- (3,0) -- (1.5,2.6) -- (0,0);
\draw (0.75,0.433) -- (2.25,0.433) -- (1.5, 1.733) -- (0.75,0.433); 
\draw (0,0) -- (0.75,0.433) -- (1.5,2.6) -- (1.5, 1.733) --(3,0) -- (2.25,0.433) -- (0,0);
\end{tikzpicture}
~
\begin{tikzpicture}
\draw (0,0) -- (3,0) -- (1.5,2.6) -- (0,0);
\draw (0.75,0.433) -- (2.25,0.433) -- (1.5, 1.733) -- (0.75,0.433); 
\draw (0,0) -- (0.75,0.433) -- (3,0) -- (2.25,0.433) -- (1.5,2.6) -- (1.5, 1.733) -- (0,0);
\end{tikzpicture}
\end{center}
\caption{Non-regular Decompositions}
\label{figure_non_triangulation}
\end{figure}
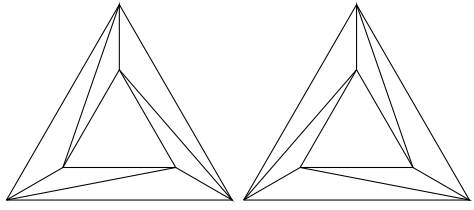

Given a pair $(\delta,A)$ regular decompositions do not always exist.  However using the theory of secondary polytopes, we can identify regular decompositions with the vertices of a certain associated polytope.  Constructions of this nature are detailed in \cite{GKZ94}.  The decompositions in Figure \ref{figure_non_triangulation} do not admit a function $\phi$ that satisfies both conditions.  Secondary polytopes can be used to show that the symmetry of these two decompositions prevents them from being regular.

\subsubsection{Indecomposable and Complete Decompositions}

A polytope $\Delta$ is called indecomposable if it contains no lattice points other than vertices.  A decomposition of a polytope $\Delta = \cup_i \Delta_i$ is called complete if each $\Delta_i$ is indecomposable.  Note a complete decomposition is always a triangulation.  In this case the $f_{\Delta_i}$ are diagonal and we may use Theorem \ref{theoremDiag}.  Complete decompositions are useful in determining whether or not $GNP(\Delta,p) = HP(\Delta)$.  

\begin{theorem}[Regular Decomposition]
Let $\cup_i \Delta_i$ be a complete regular decomposition of $\Delta$.  If each $f_{\Delta_i}$ is generically non-degenerate and ordinary for some prime $p$, then $f$ is also generically non-degenerate and ordinary for the same prime $p$.
\label{mainth}
\end{theorem}

\subsection{Newton Polygons of Subcones}
\label{section_subcones}
We now focus our attention to cones of integral polytopes.  As we saw in the previous section, $C(\Delta)$, the cone generated by $\Delta$ is useful in examining $L(f,T)$ with Dwork theory.  In particular $A_1(f)$  is deeply connected to $C(\Delta)$.  Decompositions of $\delta$ will induce a decomposition of the $C(\Delta)$ into sub-cones.  These sub-cones will not, in general, decompose $A_1(f)$ as in the case of the boundary decomposition.  We will, however, be able to replace $A_1(f)$ with a cone version that will be stringent enough to detect ordinarity.

\subsubsection{Chain Level Hodge Polygon of a Sub-cone}
Let $\Sigma$ be a cone contained in $C(\Delta)$, not necessarily open or closed.  Define a function of nonnegative integers as follows:
\begin{equation}
\label{eq64}
W(\Sigma, k)  = \# \left\{ r \in \zz^n \cap \Sigma \mid w(r) = \frac{k}{D(\Delta)} \right\}.
\end{equation}
This is the number of lattice points in the cone $\Sigma$ with weight exactly $k/D$.  Let $P(\Sigma)$ be the polygon in $\rr^2$ with vertices $(0,0)$ and 
\begin{equation}
\label{eq65}
\left( \sum_{k=0}^m W(\Sigma,k), \frac{1}{D(\Delta)} \sum_{k=0}^m kW(\Sigma,k) \right),~m=0,1,2,\ldots.
\end{equation}
For convenience, we shall call the vertex in (\ref{eq65}) the $m$\superscript{th} vertex in $P(\Sigma)$.  Note that the $m$\superscript{th} vertex may be equal to the $(m+1)$\superscript{th} vertex, because it may happen that $W(\Sigma,m)=0$.  Recall that $A_1(f)$ is the semilinear Frobenius matrix defined in (\ref{eq49}).  Recall we define $A_1(\Sigma,f)$ to be the submatrix $(a_{s,r}(f))$ with $r$ and $s$ running through the cone $\Sigma$.  For the full cone $\Sigma=C(\Delta)$, we have
$$W(C(\Delta),k) = W(k).$$
From the block form (\ref{block53}), we deduce:
\begin{prop}
The Fredholm determinant $\det(I-T A_1(\Sigma,f))$ is entire.  Its Newton polygon lies above the polygon $P(\Sigma)$.
\label{prop5.1}
\end{prop}

\subsubsection{Hasse Polynomials}
Let $P(\Sigma,x)$ be the piecewise linear function on $\rr_{\geq 0 }$ whose graph is the polygon $P(\Sigma)$.  Recall
$$f = \sum_{j=1}^J a_j V_j.$$
We may identify $f$ with its coefficients $(a_1,\ldots,a_J)$.  By the block form (\ref{block53}), we can write
\begin{equation}
\det(I-TA_1(\Sigma,f)) = \sum_{k=0}^{\infty} p^{P(\Sigma,k)} G(\Sigma,f,k)T^k,
\end{equation}
where $G(\Sigma,f,k)$ is a power series in the $a_j$ with $p$-adic integral coefficients.  The reduction
\begin{equation}
\HH(\Sigma,f,k) \equiv G(\Sigma,f,k) \modulo{\pi}
\end{equation}
is a polynomial in the coefficients $a_j$ of $f$ defined over the finite prime field $\ff_p$.  This polynomial is called the $k$\superscript{th}  Hasse polynomial of the pair $(\Sigma,f)$.

For a given pair $(\Sigma,f)$, the Newton polygon of $\det(I-TA_1(\Sigma,f))$ coincides with its lower bound $P(\Sigma)$ at the $m$\superscript{th} vertex
$$\left( \sum_{i=0}^{m} W(\Sigma,i), \frac{1}{D} \sum_{i=0}^m i W(\Sigma,i) \right)$$
if and only if the Hasse polynomial $\HH(\Sigma,f,t)$ does not vanish for 
$$k = \sum_{i=0}^m W(\Sigma,i)$$
at the point $(a_1,\ldots,a_{J})$.  To show that the Newton polygon of $\det(I-T A_1(\Sigma,f))$ coincides generically with its lower bound at the $m$\superscript{th} vertex, we need to show that the Hasse polynomial $\HH(\Sigma,f,k)$ is not identically zero for $k=\sum_{i=0}^m W(\Sigma,i)$. 

Using the Hasse polynomial, one may determine the generic Newton polygon.  However, it is very difficult to compute.  The Hasse polynomial for single variable polynomials was determined in \cite{BF07} by Blache and F\'erard.  This was done using a technique developed by Zhu in \cite{ZHU03} that is similar to the maximizing function in the next section.  A year later Liu and Chuanze generalized this result to include single variable Laurent polynomials in \cite{LIUCHU08}.


\subsection{Maximizing Functions}


Recall that $\delta$ is the unique face of $\Delta$ away from the origin.  Let $A= \{V_1,\ldots, V_J\}$ be a subset of $\delta \cap \zz^n$ which contains the vertices of $\delta$.  Let $T=\{\delta_1,\ldots, \delta_h\}$ be a decomposition of $(\delta,A)$ with an associated function $\phi:\delta \rightarrow \rr$ that is piecewise linear on each $\delta_i$ (not necessarily concave).  Therefore, for each $ r \in C(\Delta) \setminus \{0\}$ we have $r/w(r) \in \delta$.  We may naturally extend the domain of $\phi$ to $C(\Delta)$ by letting
$$\phi(r) = w(r)\phi(\frac{r}{w(r)}).$$  
If the decomposition is regular then the domains of linearity are exactly the subcones $\Sigma_i=C(\delta_i)$.  This extension of $\phi$ can be thought of as a generalization of Wan's ``priority variables" in \cite{WAN04}.  We can recover many of his original decomposition theorems through judicious construction of an appropriate $\phi$ function.  More on this is detailed at the end of this section

\begin{definition}
For $r \in C(\Delta)$ we define
	$$m(\phi,A;r) = \sup\{ \sum_{j=1}^J u_j \phi(V_j) \mid \sum_{j=1}^J u_j V_j = r, u_j \geq 0 \}.$$
	If $ r \in \rr^n$ but $r \notin C(\Delta)$, we define $m(\phi,A;r)=0$.  If for all $ r \in C(\Delta)$ we have 
	$$m(\phi,A;r) = \inf\{ \sum_{j=1}^J u_j \phi(V_j) \mid \sum_{j=1}^J u_j V_j = r, u_j \geq 0 \},$$
	we say that $\phi$ is homogeneous with respect to $A$.
	\label{def5.2}
\end{definition}

If for each $r \in A$ we set $\phi(r) = 1$ then $m(\phi,A;r)$ is the standard weight function $w(r)$ and is automatically homogeneous.  This generalized form adjusts the contributions each $V_j$ to the total weight.  Homogeneity plays a key role in decomposition theory.  

\subsubsection{An Example and Non-example of Homogeneity}
Let $\Delta$ be the convex polytope in $\rr^3$ spanned by the origin and the vectors $\smallpmatrix{2 \\ 0 \\0}$, $\smallpmatrix{0 \\ 2 \\0}$ and $\smallpmatrix{0\\0\\2}$.  Let $\delta$ be the unique codimension 1 face not containing the origin spanned by $\smallpmatrix{2 \\ 0 \\0}$, $\smallpmatrix{0 \\ 2 \\0}$ and $\smallpmatrix{0\\0\\2}$.

Let $\phi(x)=1$ for each $r=(r_1,r_2,r_3) \in \delta$ and let $A$ be the vertices of $\delta$.  Hence $m(\phi,A;r)$ is the standard weight function for $\Delta$: $w(x) = \frac{1}{2}(r_1+r_2+r_3)$.

Now, let $A =\{ \smallpmatrix{2\\0\\0},\smallpmatrix{0\\2\\0},\smallpmatrix{0\\0\\2},\smallpmatrix{ 0 \\ 1 \\ 1}\}$.  Let $\phi'(\smallpmatrix{2 \\ 0 \\0}) = \phi'(\smallpmatrix{0 \\ 1 \\1})=1 $ and $\phi'(r) =0$ for every other $r \in A$.   Then 
$$m(\phi',A;\smallpmatrix{0\\2\\2}) = 2$$ 
since $\smallpmatrix{0\\2\\2} = 2 \smallpmatrix{0\\1\\1}$ and $2\phi'(\smallpmatrix{0\\2\\2})=2$.
However we may also write,
$\smallpmatrix{0\\2\\2} = \smallpmatrix{ 0\\ 2 \\0} + \smallpmatrix{0 \\0  \\2}$
and
$\phi(\smallpmatrix{ 0\\ 2 \\0}) + \phi(\smallpmatrix{0 \\0  \\2})=0$.
  Therefore $\phi'$ is not homogeneous with respect to $A$.  Note that neither $\phi$ is not concave but $\phi'$ is.  We can make $\phi$ concave by changing the set $A$ to just the vertices.

\subsection{Degree Polygon}

For a non-negative integer $k$ and subcone $ \Sigma$ of $C(\Delta)$, we define 
\begin{equation}
Q(\Sigma, \phi, A;k) = (p-1) \sum_{w(r) \leq k/D, r \in \zz^n\cap \Sigma} m(\phi, A;r)
\label{eq67}
\end{equation}
where the intersection $\zz^n \cap \Sigma$ is simply the set of lattice points in the cone $\Sigma$.  Let $Q(\Sigma,\phi,A)$ be the graph in $\rr^2$ of the piecewise linear functions passing through the vertices $(0,0)$ and 
\begin{equation}
\left(\sum_{k=0}^m W(\Sigma,k), Q(\Sigma, \phi ,A;m) \right), m=0,1,2,\ldots.
\label{eq68}
\end{equation}
We shall call the vertex in (\ref{eq68}) the $m$\superscript{th} vertex in $Q(\Sigma, \phi,A)$.  Note that the coordinates in (\ref{eq68}) are always non-negative.  In the special case that $\Sigma = C(\Delta)$ we simply write
$$Q(C(\Delta), \phi,A) = Q(\phi,A).$$

Recall $F_r$ is the polynomial in the variables $a_j$ defined in (\ref{eq45}):
$$F_r(f) = \sum_u \left( \prod_{j=1}^J \lambda_{u_j} a_j^{u_j}  \right) \pi^{u_1+\ldots+u_J},$$
where the outer sum is over all solutions of the linear system
$$\sum_{j=1}^J u_j V_j =r , u_j \geq 0, u_j \mathrm{~integral}.$$

\begin{definition}
For any polynomial 
$$F(a_1,\ldots,a_J) = \sum_{u=(u_1,\ldots,u_j)} \lambda_u \prod_{j=1}^J a_j^{u_j},$$
 define $d(\phi,F)$ to be the total degree of $F$, weighted by the function $\phi$.  That is,
$$d(\phi,F) = \max \{ \sum_{j=1}^J u_j \phi(V_j) \},$$
where the maximizing function is taken over all $u=(u_1,\ldots,u_j)$ such that $\lambda_u \neq 0$.  If $F=0$, we define $d(\phi,F)=-\infty$.  The number $d(\phi,F)$ is called the $\phi$-degree of $F$.  
\end{definition}

In the case that $\phi$ is homogeneous with respect to $A$, $F_r$ is $\phi$-homogeneous and the above inequality becomes an equality if and only if $F_r$ is non-zero.

\begin{definition}
Let $d(\phi, \HH(\Sigma, f,k))$ denote the $\phi$-degree of the $k$\superscript{th} Hasse polynomial $\HH (\Sigma,f,k)$.  We define the $\phi$-degree polygon of $\det(I-TA_1(\Sigma,f))$ to be the graph in $\rr^2$ of the piecewise linear function with vertices $(0,0)$ and 
$$\left( \sum_{k=0}^m W(\Sigma,k), \max \{ 0, d(\phi, \HH (\Sigma,f, \sum_{k=0}^m W(\Sigma,k)))\} \right), m=0,1,\ldots,\infty.$$
Note that we do not claim that the $\phi$-degree polygon is convex.
\end{definition}

If $r$ and $r'$ are two lattice points in $\zz^n$, Definition \ref{def5.2} implies the inequality
\begin{equation}
m(\phi,A;r) + m(\phi,A;r') \leq m(\phi,A;r+r').
\label{eq70}
\end{equation}
The equality always holds if $\phi$ is homogeneous with respect to $A$ and $r,r' \in C(\Delta)$.  Furthermore, if $c$ is non-negative, then
$$m(\phi,A;cr) = cm(\phi,A;r).$$
Let $m$ be a non-negative integer.  By the block form of the matrix $A_1(\Sigma,f)$, the determinant expansion of a matrix and the development in this section, we deduce
$$d \left( \phi, \HH (\Sigma,f, \sum_{k=0}^m W(\Sigma,k)) \right)$$
$$\begin{array}{rl}
\leq &\displaystyle{\max_{\psi}} \sum_{w(r) \leq m/D,r\in \zz^n \cap \Sigma} d(\phi,F_{pr-\psi(r)})\\
\leq &\displaystyle{\max_{\psi}} \sum_{w(r) \leq m/D,r\in \zz^n \cap \Sigma} m(\phi,A;pr-\psi(r))\\
\leq &\displaystyle{\max_{\psi}} \sum_{w(r) \leq m/D,r\in \zz^n \cap \Sigma} ( (pm(\phi,A;r) -  m(\phi,A;\psi(r)) )\\
\leq& (p-1) \sum_{w(r) \leq m/D,r\in \zz^n \cap \Sigma} m(\phi,A;r)\\
=&Q(\Sigma, \phi, A;m)
\end{array}$$
where $\psi$ runs through the permutations on $\sum_{k=0}^m W(\Sigma,k)$ letters.  If $\phi$ is homogeneous with respect to $A$ and if the Hasse polynomial
$$\HH (\Sigma,f, \sum_{k=0}^m W(\Sigma,k))$$
is not the zero polynomial, then the Hasse polynomial is $\phi$-homogeneous and we have the equality
$$d \left( \phi, \HH (\Sigma,f, \sum_{k=0}^m W(\Sigma,k)) \right)=Q(\Sigma, \phi, A;m).$$

\begin{prop}
The $\phi$-degree polygon of $\det(I-T A_1(\Sigma,f))$ lies below the polygon $Q(\Sigma, \phi, A)$.
\end{prop}
If the $\phi$-degree polygon coincides with its upper bound $Q(\Sigma, \phi,A)$ at the $m$\superscript{th} vertex, then the polynomial $\HH (\Sigma,f, \sum_{i=0}^m W(\Sigma,i))$ is not zero since its $\phi$-degree is equal to $Q(\Sigma, \phi, A;\sum_{i=0}^m W(i))$ which is non-negative and hence not $-\infty$.  The converse is also true if $\phi$ is homogeneous with respect to $A$.  From this we obtain the following proposition.
\begin{prop}
If the $\phi$-degree polygon of $\det(I - TA_1(\Sigma,f))$ coincides with $Q(\Sigma, \phi, A)$ at the $m$-th vertex, then the Newton polygon of $\det(I-T A_1(\Sigma,f))$ coincides generically with $P(\Sigma)$ at the $m$-th vertex.  If $\phi$ is homogeneous with respect to $A$, then the converse is also true.
\label{prop5.6}
\end{prop}
This property shows that the $\phi$-degree polygon is finer than the generic Newton polygon.  When the $\phi$-degree polygon of $\det(I - TA_1(\Sigma,f))$ coincides with $Q(\Sigma, \phi, A)$ we call this ordinarity on the degree polygon level.


\section{Decomposition Theorems}
\label{chapter_decomposition}


\subsection{Regular Decompositions}
Consider a $T=\{\delta_1,\ldots,\delta_h\}$ be regular decomposition of $\delta$,
with associated concave function $\phi$ and set of vertices $A$.  For each $1\leq i \leq h$, let $\Sigma_i = C(\delta_i)$, the cone generated by $\delta_i$ and the origin.  


To prove Theorem \ref{mainth} we will first prove a version for degree polygons.   This closed regular decomposition will itself require several steps.  We first state the theorem, then use it to prove Theorem \ref{mainth}.
\begin{theorem}
\label{th6.1}
For $m \in \zz_{\geq 0}$, the $\phi$-degree polygon of $\det(I-T A_1(f))$ coincides with its upper bound $Q(\phi,A)$ at the $m$\superscript{th} vertex if and only if  for each $1 \leq i \leq h$, the $\phi|_{\delta_i}$-degree polygon of $\det(I-T A_1(f_{\Sigma_i}))$ defined with respect to $P(\Sigma_i)$ coincides with its upper bound $Q(\phi|_{\delta_i},\delta_i)$ at the $m$\superscript{th} vertex.  Here $\phi|_{\delta_i}$ denotes the restricted function where $\phi|_{\delta_i}=\phi$ at points in $\delta_i$ and vanishes elsewhere.
\label{thclosedregular}
\end{theorem}

If $|A|=n$, the minimal possible value, then no additional decomposition is possible and there is nothing to prove.  Suppose $|A|>n$.  Under the assumptions of Theorem \ref{mainth} 
the Newton polygon of $\det(I-TA_1(f_{\Sigma_i}))$ coincides generically with $P(\Sigma_i)$.  Since $\delta_i$ is indecomposable, $\phi|_{\delta_i}$ is homogeneous with respect to $\delta_i$.  By Proposition \ref{prop5.6}, the $\phi|_{\delta_i}$-degree polygon of $\det(I-TA_1(f_{\Sigma_i}))$ defined with respect to $P(\Sigma_i, f_{\Sigma_i})$ coincides with its upper bound $Q(\phi_{\delta_i},\delta_i)$.  By Theorem \ref{th6.1}, we deduce that the $\phi$-degree polygon of $\det(I-TA_1(f))$ coincides with its upper bound $Q(\phi,A)$.  Applying Proposition \ref{prop5.6} again, we conclude that $f$ is generically ordinary.  Theorem \ref{mainth} is proved.

Note that completeness of the decomposition is necessary to ensure homogeneity.  Without homogeneity we could not use Proposition \ref{prop5.6} to move from the chain level to the degree polygon level.

Before proving Theorem \ref{th6.1}, we need to have a better understanding of the maximizing function that is used to define the degree polygon.  Recall that for $r \in C(\Delta)$, we defined
$$m(\phi,A;r) = \sup\{ \sum_{j=1}^J u_j \phi(V_j) \mid \sum_{j=1}^J u_j V_j = r, u_j \geq 0 \}.$$
If $r \in \rr^n$ but $r \notin C(\Delta)$, then $m(\phi,A;r)=0$.
\begin{lemma}
Let $T=\{\delta_1,\ldots,\delta_h\}$ be a regular decomposition of $(\delta,A)$ with associated concave function $\phi:\delta \rightarrow \rr$.  Let $(u_1,\ldots,u_J)$ be a rational solution of the linear equation
$$\sum_{j=1}^J u_jV_j=r, u_j \geq 0.$$
Suppose that $u_{j_1},\ldots,u_{j_k}$ are its non-zero coordinates.  If $r \in \Sigma_i$, then 
$$m(\phi,A;r)=\sum_{j=1}^J u_j \phi(V_j)$$
 if and only if $\Sigma_i$ contains all the lattice points $V_{j_1},\ldots,V_{j_k}$.
\label{lemma6.2}
\end{lemma}
\it Proof. \rm  If $u_{j_1},\ldots,u_{j_k}$ are the non-zero coordinates of $r$ and some of the $V_{j_1},\ldots, V_{j_k}$ are all contained in $\Sigma_i$ then we may express $r$ as the sum $r_1+r_2$ where $r_1\in \Sigma_i$ and $r_2 \notin \Sigma_i$.  Note that the choices for $r_1$ and $r_2$ are not necessarily unique.  Since the decomposition is regular and $r_1$ and $r_2$ lie in different domains of linearity, we have
$$\phi(r) = \phi(r_1+r_2) > \phi(r_1) + \phi(r_2).$$
Therefore $\phi$ is not maximized and $m(\phi, A;r)>\phi(r_1) + \phi(r_2)$.  A similar exercise in the definitions will show that if $\phi$ is maximized then  $V_{j_1},\ldots,V_{j_k}$ are in $\Sigma_i$.  The proof is complete.

Using this lemma we obtain another useful lemma.

\begin{lemma}
Let $r_1$ and $r_2$ be two rational points in the cone $C(\Delta)$.  Then
\begin{equation}
m(\phi,A;r_1+r_2) \geq m(\phi,A,r_1) + m(\phi,A,r_2).
\label{lemma_meq}
\end{equation}
If equality holds, then both $r_1$ and $r_2$ lie on $\Sigma_i$ for some $i$.
\label{lemma6.3}
\end{lemma}
\it Proof. \rm Let 
$$r_1 = u_1 V_1 + \ldots +u_J V_J,~ \sum_{j=1}^J u_j \phi(V_j)= m(\phi,A;r_1),$$
$$r_2 = w_1 V_1 + \ldots +w_J V_J,~ \sum_{j=1}^J w_j\phi(V_j) = m(\phi,A;r_2).$$
Then,
$$r_1+r_2 = (u_1+w_1)V_1 + \ldots + (u_J+w_J)V_J.$$
$$m(\phi,A;r_1)+m(\phi,A;r_2) = \sum_{j=1}^J (u_j+w_j) \leq m(\phi,A;r_1+r_2).$$
If equality holds in (\ref{lemma_meq}), then by Lemma \ref{lemma6.2}  the $V_j$ with nonzero coefficients all appear in $\Sigma_i$.  Since the coefficients $u_i$ and $w_i$ are non-negative, the set of lattice points with nonzero coefficients of $r_1+r_2$, contain the set of nonzero lattice points of both $r_1$ and $r_2$.  The lemma is proved.

The next lemma shows that certain leading terms of $F_r(f)$ and $F_r(f_{\Sigma})$ are identical.
\begin{lemma}
Let $r \in \Sigma_i$ for some $1 \leq i \leq h$ and $f=\sum_{j=1}^J a_j x^{V_j}$.  Then 
$$d(\phi, F_r(f) - F_r(f_{\Sigma_i})) < m(\phi, A;r).$$
That is, the $\phi$-degree of $F_r(f) - F_r(F_{\Sigma_i})$ is strictly smaller than the expected maximum value $m(\phi,A;r)$.
\label{lemma6.4}
\end{lemma}

\noindent\it Proof. \rm
Let $u_1,\ldots,u_J$ be non negative integers satisfying
$$r=\sum_{j=1}^J u_j V_j, \sum_{j=1}^J u_j \phi(V_j)=m(\phi,A;r).$$
Let $u_{j_1}, \ldots, u_{j_k}$ be the non-zero terms among the $u_j$'s.  Lemma \ref{lemma6.2} show that the cone $\Sigma_i$ contains all the lattice points $V_{j_1},\ldots, V_{j_k}$.  This shows that if a monomial in the $a_j$ with $\phi$-degree $m(\phi,A;r)$ appears in $F_r(f)$, then the same monomial also appears in $F_r(f_{\Sigma_i})$.  Thus $F_r(f)$ and $F_r(f_{\Sigma_i})$ have the same initial terms of $\phi$-degree $m(\phi,A;r)$.  The lemma is proved.



For $1 \leq i \leq h$, let $S_i^o$ be the set of relatively open faces in $\delta_i$, including the empty set.  This is also called the open boundary decomposition of $\delta_i$.  The union.
$$\delta_i = \bigcup_{\sigma \in S_i^o} \sigma$$
is a disjoint union.  Each $S_i^o$ contains exactly one $(n-1)$-dimensional face, namely the interior of $\delta_i$.  The $0$-dimensional elements in $S_i^o$ are simply the vertices in $\delta_i$.  For an element $\sigma \in S_i^o$, let $\Sigma^o$ denote the open cone generated by $\sigma$ and the origin.  The origin itself is not included unless $\sigma$ is the empty set.  Since $\sigma \in S_i^o$, the open cone $\Sigma^o$ is a subcone of the closed cone $\Sigma_i$ and we have
$$\dim \Sigma^o = \dim \sigma+1.$$
The union
$$\Sigma_i = \bigcup_{\sigma \in S_i^o} \Sigma^o$$
is a disjoint union, called the boundary decomposition of $\Sigma_i$.  Let $\Sigma$ be the topological closure of $\Sigma^o$.  To prove the closed regular decomposition we first need to prove an open version:
\begin{theorem}
For $m \in \zz_{\geq0}$ the $\phi$-degree polygon of $\det(I-TA_1(f))$ coincides with its upper bound $Q(\phi,A)$ at the $m$\superscript{th} vertex if and only if for each $1 \leq i \leq h$ and each $\sigma \in S_i^o$, the $\phi|_{\sigma}$-degree polygon of $ \det(I-T A_1(\Sigma^o,f_{\Sigma_i}))$ coincides with its upper bound $Q(\Sigma^o, \phi|_{\Sigma_i}, \Sigma_i)$ at the $m$\superscript{th} vertex.
\label{opendecomptheorem}
\end{theorem}
This is another instance of a property holding for the components if and only if that property hold for the whole.

\noindent
\it Proof. \rm
Let 
$$S(\phi,A) = \bigcup_{i=1}^h S_i^o$$
which includes the empty set as an element.  Denote the number of elements in this set by $g+1$.  Fix an ordering of $S(\phi,A)$ by
$$S(\phi,A) = \{\sigma_0, \sigma_1,\ldots,\sigma_g\}$$
such that 
$$\dim(\sigma_j) \leq \dim(\sigma_{j+1}), 0 \leq j \leq g+1.$$
In particular, $\sigma_0$ is the empty set.  Let $\Sigma_j^o$ be the relatively open cone generated by $\sigma_j$ and the origin.  In particular, $\Sigma_0^o$ consists of the origin.  It is clear that we have 
$$\dim(\Sigma_j^o) = \dim(\sigma_j)+1$$
and thus
$$\dim(\Sigma_j^o) \leq \dim(\Sigma_{j+1}^o), 0 \leq j \leq g-1.$$
Let $$C(\phi,T) = \{\Sigma_0^o,\ldots,\Sigma_g^o\}.$$
Thus the full cone $C(\Delta)$ is the disjoint union of the relatively open cones in $C(\phi,A)$:
$$C(\Delta) = \cup_{j=0}^g \Sigma_j^o.$$
Let $0 \leq j_1 \leq j_2 \leq g$ and
$$s \in \Sigma_{j_1}^o, r \in \Sigma_{j_2}^o.$$
In particular, $r \neq 0$ since the origin is only contained in $\Sigma_0^o$.  We claim that for each $1 \leq i \leq h$, the closed convex cone $\Sigma_i$ cannot contain both $r$ and $ps-r$.  Otherwise, suppose that both $r$ and $ps-r$ are contained in $\Sigma_i$ for some $1 \leq i \leq h$.  Observe $ps=r+(ps-r)$.  It follows that $r$, $s$ and $ps-r$ are all contained in $\Sigma_i$.  In particular, both $\Sigma_{j_1}^o$ and $\Sigma_{j_2}^o$ are subcones of $\Sigma_i$.  Let $ps-r \in \Sigma_{j_3}^o$ (a subcone of $\Sigma_i$).  Then the equation $ps=r+(ps-r)$ shows that $s$ is in the interior of the cone generated by $\Sigma_{j_2}^o$ and $\Sigma_{j_3}^o$.  This implies that 
$$\dim(\Sigma_{j_1}^o) \geq \dim(\Sigma_{j_2}^o),$$
with equality holding if and only if $ps-r$ lies in the topological closure of $\Sigma_{j_2}^o$.  Our ordering assumption shows that this is indeed an equality.  Thus $ps-r$ is indeed in the topological closure of $\Sigma_{j_2}^o$.  We conclude from $ps = r +(ps-r)$ that $s$ is in $\Sigma_{j_2}^o$.  This shows that $\Sigma_{j_1}^o$ and $\Sigma_{j_2}^o$ are not disjoint, a contradiction.  The claim is proved.  This claim together with Lemma \ref{lemma6.3} shows that
\begin{equation}
\begin{array}{rcccl}
d(\phi,a_{s,r}(f)) & = & d(\phi,F_{ps-r}(f))& &\\
&\leq & m(\phi,A;ps-r)& <& pm(\phi, A;s) - m(\phi,A;r).
\end{array}
\label{eq86}
\end{equation}

Let $B_{j_1 j_2} ( 0 \leq j_1, j_2 \leq g)$ be the nuclear submatrix of $A_1(f)$ consisting of all $(a_{s,r}(f))$ with $s \in \Sigma_{j_1}^o$ and $r \in \Sigma_{j_2}^o$.  For $0 \leq j \leq g$, the Newton polygon of the entire function $\det(I-tB_{jj})$ lies above $P(\Sigma_j^o)$.  Furthermore, under a permutation of orthonormal basis we see that $A_1(f)$ is similar to the matrix:
\begin{equation}
\begin{pmatrix} 
B_{00} & B_{01} & \ldots & B_{0g} \\
B_{10} & B_{11} & \ldots & B_{1g} \\
\vdots & \vdots & \ddots & \vdots \\
B_{g0} & B_{g1}& \ldots & B_{gg}
 \end{pmatrix} .
 \label{eq87}
 \end{equation}
 
 If $a_{s,r}$ is an element in $B_{j_1 j_2}$, then (\ref{eq86}) shows that the $\phi$-degree $d(\phi,a_{s,r}(f))$ of the polynomial $a_{s,r}(f)$ is strictly smaller than the expected maximum value $pm(\phi,A;s)-m(\phi,A;r)$.  This means that the block form for $A_1(f)$ is in some sense, lower triangular with respect to the $\phi$-degree.  By induction we deduce that
 \begin{equation}
 \det(I-TA_1(f)) = \prod_{j=0}^g \det(I-TB_{jj}) + \sum_{k=0}^{\infty} p^{P(\Delta,k)} G(f,k)T^k,
 \label{eq88}
 \end{equation}
 where $G(f,k)$ is a power series in the $a_j$ with $p$-adic integral coefficients.  The reduction $G(f,k) \modulo{\pi}$ is a polynomial over $\ff_p$ whose $\phi$-degree is strictly smaller than the upper bound $Q(\phi,A;k)$ (see the notation in (\ref{eq67}) and (\ref{eq68})).  Thus the $\phi$-degree polygon of $\det(I-TA_1(f))$ coincides with $Q(\phi,A)$ at the $m$\superscript{th} vertex if and only if the $\phi$-degree polygon of the first term on the right side of $(\ref{eq88})$ coincides with $Q(\phi,A)$ at the $m$\superscript{th} vertex.  One further shows that the latter is true if and only if the $\phi$-degree polygon of
 $$\det(I-TA_1(\Sigma_j^o,f)) = \det(I-TB_{jj})$$
 defined with respect to $P(\Sigma_j^o)$ coincides with $Q(\Sigma_j^o,\phi,A)$ at the $m$\superscript{th} vertex for all $0\leq j\leq g$.  The matrix $A_1(\Sigma_j^o,f)$ is however different from the desired matrix $A_1(\Sigma_j^o,f_{\Sigma_i})$, where $\Sigma_j^o \subset \Sigma_i$.  But the $\phi$-degree polygon of $\det(I-TA_1(\Sigma_j^o),f)$ and the $\phi|_{\delta_i}$-degree polygon of $\det(I-TA_1(\Sigma_j^o,f_{\Sigma_i}))$ have the same upper bound
 \begin{equation}
 Q(\Sigma_j^o,\phi,A) = Q(\Sigma_j^o, \phi|_{\delta_i},\delta_i).
 \label{eq89}
 \end{equation}
The last equality can be proved from our definitions in (\ref{eq67}) and (\ref{eq68}).

To finish the proof, we need to show that if we replace the matrix $A_1(\Sigma_j^o,f)$ defined in terms of $f$ by the matrix $A_1(\Sigma_j^o,f_{\Sigma_i})$ defined in terms of $f_{\Sigma_i}$, we will not change the property of the coincidence of the degree polygon with its upper bound.  Let $r,s \in \Sigma_i$.  If $ps-r$ also belongs to $\Sigma_i$, then Lemma \ref{lemma6.4} shows that we can replace $F_{ps-r}(f)$ by $F_{ps-r}(f_{\Sigma_i})$.  If $ps-r$ does not belong to $\Sigma_i$, then Lemma \ref{lemma6.3} shows that the $\phi$-degree polygon of $F_{ps-r}(f)$ is strictly smaller than the expected maximum value $pm(\phi,A;s)-m(\phi,A;r)$, while $F_{ps-r}(f_{\Sigma_i})=0$.  In this case, we can also replace $F_{ps-r}(f)$ by $F_{ps-r}(f_{\Sigma_i})$.  Hence Theorem \ref{opendecomptheorem}, the open regular decomposition theorem, is proved.

To prove the closed regular decomposition theorem, it suffices to combine the above open regular decomposition theorem and  apply Theorem \ref{th5.1}, the boundary decomposition theorem, to each $\Sigma_i$.

Theorem \ref{th5.1} shows that the $\phi|_{\delta_i}$-degree polygon of $\det(I-TA_1(f_{\Sigma_i}))$ coincides with its upper bound $Q(\Sigma_i, \phi|_{\delta_i}, \delta_i)$ if and only if the $\phi|_{\delta_i}$-degree polygon of $\det(I-TA_1(\Sigma_j^o,f_{\Sigma_i}))$ coincides with its upper bound $Q(\Sigma_j^o,\phi|_{\delta_i},\delta_i)$ for all $j$ with $\sigma_j \in S_i^o$.  The proof of Theorem \ref{th6.1} is complete.

\subsection{Other Decomposition Theorems}

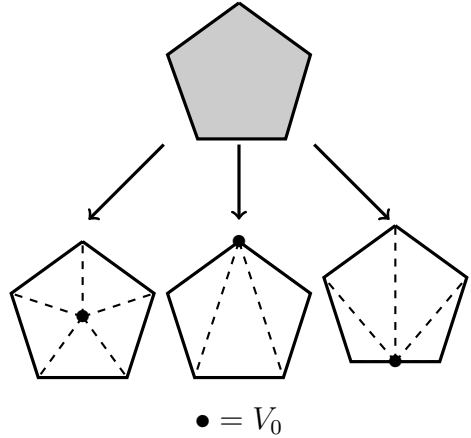
\begin{figure}[t]
\begin{center}
\begin{tikzpicture}
\fill[black!20!white] (0, 1.0)--(-0.95, 0.31)--(-0.55, -0.81)--(0.62, -0.81)--(0.95, 0.31) -- (0,1);
\draw[very thick] (0, 1.0)--(-0.95, 0.31)--(-0.55, -0.81)--(0.62, -0.81)--(0.95, 0.31) -- (0,1);
\end{tikzpicture}
\\
\begin{tikzpicture}
\draw[-to,very thick] (-1,0) -- (-2,-1);
\draw[-to,very thick] (0,0) -- (0,-1);
\draw[-to,very thick] (1,0) -- (2,-1);
\end{tikzpicture}
\\
\begin{tikzpicture}
x=(1,1);
\draw[very thick] (0, 1.0)--(-0.95, 0.31)--(-0.59, -0.81)--(0.59, -0.81)--(0.95, 0.31) -- (0,1);
\draw[dashed,thick] (0,0) -- (0,1);
\draw[dashed,thick] (0,0) -- (-0.95, 0.31);
\draw[dashed,thick] (0,0) -- (-0.59, -0.81);
\draw[dashed,thick] (0,0) -- (0.59, -0.81);
\draw[dashed,thick] (0,0) -- (0.95, 0.31);
\draw (0,0) node {$\bullet$};
\end{tikzpicture}
~
\begin{tikzpicture}
\draw[very thick] (0, 1.0)--(-0.95, 0.31)--(-0.59, -0.81)--(0.59, -0.81)--(0.95, 0.31) -- (0,1);
\draw[dashed,thick] (0,1) -- (0,1);
\draw[dashed,thick] (0,1) -- (-0.95, 0.31);
\draw[dashed,thick] (0,1) -- (-0.59, -0.81);
\draw[dashed,thick] (0,1) -- (0.59, -0.81);
\draw[dashed,thick] (0,1) -- (0.95, 0.31);
\draw (0,1) node {$\bullet$};
\end{tikzpicture}
~
\begin{tikzpicture}
\draw[very thick] (0, 1.0)--(-0.95, 0.31)--(-0.59, -0.81)--(0.59, -0.81)--(0.95, 0.31) -- (0,1);
\draw[dashed,thick] (0,-0.81) -- (0,1);
\draw[dashed,thick] (0,-0.81) -- (-0.95, 0.31);
\draw[dashed,thick] (0,-0.81) -- (-0.59, -0.81);
\draw[dashed,thick] (0,-0.81) -- (0.59, -0.81);
\draw[dashed,thick] (0,-0.81) -- (0.95, 0.31);
\draw (0,-0.81) node {$\bullet$};
\end{tikzpicture}
$$\bullet = V_0$$

\end{center}
\caption{Three Examples of the Star Decomposition}
\label{figurestardecomp}
\end{figure}

This section includes other decomposition theorems.  The star decomposition and the parallel hyperplane decomposition appear in \cite{WAN93}.  The collapsing decomposition appears in  \cite{WAN04}.  We will show how each of these decomposition theorems can be realized as a regular decomposition.  Using the facial decomposition we may assume that $\Delta$ contains a unique face $\delta$ of codimension $1$ not containing the origin.  Recall that these decompositions are decompositions of $\delta$, which induce a decomposition of $\Delta$.

\subsubsection{Star Decomposition}
Let $V_1,\ldots,V_J$ be the $J$ vertices of $\delta$.  Let $V_0$ be a lattice point in $\delta$ (possibly a vertex).  Let $\sigma_1,\ldots,\sigma_k$ be the open faces of codimension $1$ in $\delta$ that do not contain $V_0$.  Let $\delta_i$ be the closed, convex closure of $V_0$ and $\sigma_i$.  Let $\Delta_i$ be the closed convex closure of $\delta_i$ and the origin.  The decomposition
$$\Delta = \bigcup_i \Delta_i$$
is called the star decomposition of $\Delta$.
This can be made a regular triangulation in the following way.  Let $A = \{V_0\} \cup \{V_1,\ldots,V_J\}$.  Define $\psi:A \rightarrow \rr$ by setting $\psi(V_0) =1 $ and $\psi(V_i)=0$ for all vertices such that $V_0 \neq V_i$.  The function $g_{\psi}:\delta \rightarrow \rr$ defined in section \ref{chapter_poly} induces a regular decomposition $T$ of $(\delta,A)$.  

We can determine $g_{\psi}$ constructively.  For $r\in \delta_i$, define $d(\sigma_i,r)$ to be the distance from $r$ to $\sigma_i$ on $\delta$.  Then
$$g_{\psi}(r) = \frac{d(\delta_i,r)}{d(\delta_i,V_0)}.$$
The denominator $d(H_i,V_0)$ ensures continuity across the entire domain $\delta$ and normalizes the function.
Figure \ref{figurestardecomp} shows the star decomposition of a pentagon for various choices of $V_0$.



\subsubsection{Parallel Hyperplane Decomposition}
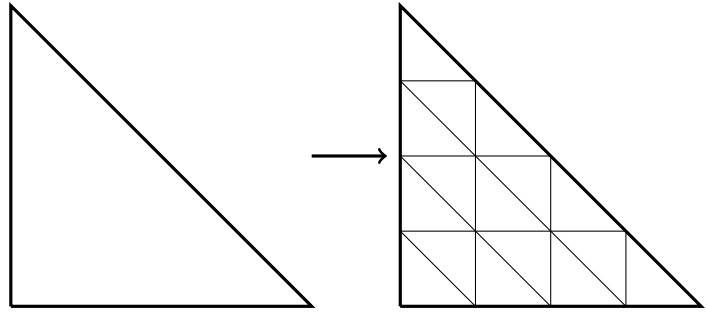
\begin{figure}[t]
\begin{center}
\begin{tikzpicture}
\draw[very thick] (0,0) -- (4,0) -- (0,4) -- (0,0);
\draw[-to, very thick] (4,2) -- (5,2);
\end{tikzpicture}
\begin{tikzpicture}
\draw[very thick] (0,0) -- (4,0) -- (0,4) -- (0,0);
\draw (0,1) -- (3,1);
\draw (0,2) -- (2,2);
\draw (0,3) -- (1,3);
\draw (1,0) -- (1,3);
\draw (2,0) -- (2,2);
\draw (3,0) -- (3,1);
\draw (1,0) -- (0,1);
\draw (2,0) -- (0,2);
\draw (3,0) -- (0,3);
\end{tikzpicture}
\end{center}
\caption{Example of the Parallel Hyperplane Decomposition}
\label{figureplanedecomp}
\end{figure}

Let $H$ be a hyperplane with codimension 1 such that the intersection of $\delta$ and $H$ is a polytope of codimension $2$ with integral vertices.  This hyperplane ``cuts" the face $\delta$ into two polytopes $\delta_1$ and $\delta_2$.  This division induces a decomposition of $\Delta$ into two polytopes $\Delta_1$ and $\Delta_2$.

To realize $\Delta$ as a regular decomposition we use a procedure similar to the case of star decompositions.  Let $V_1,\ldots,V_J$ be the vertices of $\delta$.  Let $V_{1,1}, \ldots, V_{1,k}$ be the vertices of $\delta_1$ and $V_{2,1}, \ldots, V_{2,k'}$ be the vertices of $\delta_2$.  Next, let $A = \{V_{1,1}, \ldots, V_{1,k}\} \cup \{V_{2,1}, \ldots, V_{2,k'}\}$, the union of vertices of $\delta_1$ and $\delta_2$.  Let $\phi(V)=1$ for $V \in \{V_{1,1}, \ldots, V_{1,k}\} \cap \{V_{2,1}, \ldots, V_{2,k'}\}$; these are the lattice points on $H$.  Define a constant $d_{max} = \max_{r \in \delta} (d(H,r))$.  This is the maximum distance of points on $\delta$ from $H$.  The function
$$g_{\psi}(r) = 1-\frac{d(H,r)}{d_{max}}$$
is a convex function and the decomposition $\delta = \delta_1 \cup \delta_2$ is regular.  Note that if we used $\frac{d(H,r)}{d_{max}}$ instead, this would be a convex function.  This construction is slightly different than the construction used in \cite{WAN04}.  The corresponding construction in Wan's work would have set $g_{\psi}(r) = 1$ for $r \in \delta_2$.


Parallel hyperplanes can be used form a parallel hyperplane decomposition. 
This decomposition is very useful in situations where several parallel hyperplanes apply to the same $\delta$.  Figure \ref{figureplanedecomp} shows how a typical polytope can be subdivided using three sets of parallel hyperplanes.

\subsubsection{Collapsing Decomposition}
Let $A= \{V_1,\ldots, V_J\}$ be the set of $J$ fixed lattice points in $\delta$, which include the vertices of $\delta$.  Choose an element of $A$ which is a vertex of $\delta$, say $V_1$.  Let $A_1 = A \setminus \{V_1\}$, the complement of $V_1$ in $A$.  Let $\delta_1$ be the convex polytope generated by the lattice points in $A_1$.  This is a subset of $\delta$.  Let $\delta_1'$ be the topological closure of $\delta-\delta_1$.  This is not a convex polyhedron in general.  The intersection $\delta_1 \cap \delta_1'$ consists of finitely many different codimension 2 faces $\{\sigma_2,\ldots,\sigma_h\}$ of $\delta_1$.  Let $\delta_i (2 \leq i \leq h)$ be the convex closure of $\sigma_i$ and $V_1$.  Then, each $\delta_i$ is $(n-1)$-dimensional.  The collapsing decomposition is defined to be 
$$\delta=\delta_1 \cup \ldots \cup \delta_h.$$
Let $A_i (2 \leq i \leq h)$ be the intersection of $A \cap \delta_i$.  Then, each $V_j$ lies in at least one (possibly more) of the subsets $A_i$ of $A$.  We also have a collapsing decomposition of the lattice points $A$ with respect to $V_1$: 
$$A= \bigcup_{i=1}^h A_i.$$
The collapsing decomposition can be made regular by setting $\psi(V_0)=0$ and $\psi(V_i)=1$ for all other lattice points.  The function $g_{\psi}$ can be computed explicitly:
$$g_{\psi}(r) = 
\begin{cases}
1 & \mathrm{~if~} r \in \delta_1',\\
1-\frac{d(\omega_i,r)}{d(\omega_i,V_1)}. & \mathrm{~if~} r \in \delta_i \setminus \delta_1 \mathrm{~and~} 2 \leq i \leq h.
\end{cases}$$

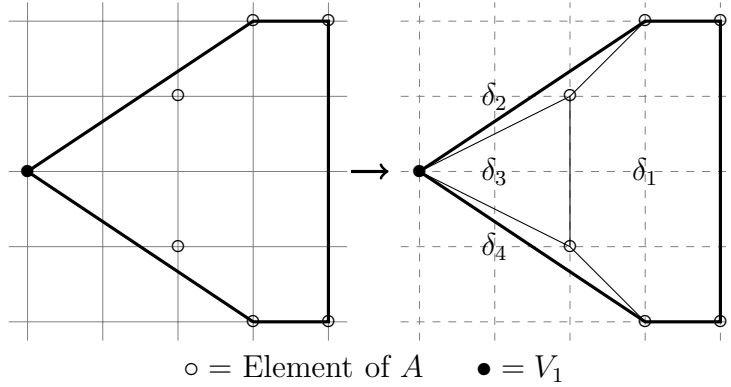
\begin{figure}[t]
\begin{center}
\begin{tikzpicture}
\draw[step= 1cm,gray,very thin] (-0.25,-0.25) grid (4.25,4.25);
\draw[very thick] (0,2) -- (3,4) -- (4,4) -- (4,0) -- (3,0) -- (0,2);
\draw (0,2) node {$\bullet$};
\draw (3,0) node {$\circ$};
\draw (2,1) node {$\circ$};
\draw (2,3) node {$\circ$};
\draw (3,4) node {$\circ$};
\draw (4,4) node {$\circ$};
\draw (4,0) node {$\circ$};
\draw[-to, very thick] (4.30,2) -- (4.8,2);
\end{tikzpicture}
\begin{tikzpicture}
\draw[step= 1cm,gray,very thin,dashed] (-0.25,-0.25) grid (4.25,4.25);
\draw[very thick] (0,2) -- (3,4) -- (4,4) -- (4,0) -- (3,0) -- (0,2);
\draw (0,2) node {$\bullet$};
\draw (3,0) node {$\circ$};
\draw (2,1) node {$\circ$};
\draw (2,3) node {$\circ$};
\draw (3,4) node {$\circ$};
\draw (4,4) node {$\circ$};
\draw (4,0) node {$\circ$};
\draw (0,2) -- (2,3);
\draw (0,2) -- (2,1);
\draw (3,4) -- (2,3) -- (2,1) -- (3,0);
\draw (3,2) node {$\delta_1$};
\draw (1,3) node {$\delta_2$};
\draw (1,2) node {$\delta_3$};
\draw (1,1) node {$\delta_4$};
\end{tikzpicture}

\begin{tabular}{l c r}
$\circ=$   Element of $A$&~~~~~ & $\bullet = V_1$
\end{tabular}

\end{center}
\caption{Example of the Collapsing Decomposition}
\label{figurecollapsingdecomp}
\end{figure}

Figure \ref{figurecollapsingdecomp} illustrates a decomposition into four pieces by collapsing at one point.

\subsubsection{Decomposition Theorems}
Each of the three decompositions in this section lead to a decomposition of ordinarity on the degree polygon level as in the closed collapsing decomposition theorem for degree polygons.  
\begin{theorem}
Let the $T=\{\delta_1,\ldots,\delta_h\}$ 
be a star, parallel hyperplane, or collapsing decomposition.  Then there is a function $\phi:\delta \rightarrow \rr$ and set of lattice points $A$ that makes $T$ a regular decomposition.  For $m \in \zz_{\geq 0}$, the $\phi$-degree polygon of $\det(I-T A_1(f))$ coincides with its upper bound $Q(\phi,A)$ at the $m$\superscript{th} vertex if and only if  for each $1 \leq i \leq h$, the $\phi|_{\delta_i}$-degree polygon of $\det(I-T A_1(f_{\Sigma_i}))$ defined with respect to $P(\Sigma_i,f_{\Sigma_i})$ coincides with its upper bound $Q(\phi|_{\delta_i},\delta_i)$ at the $m$\superscript{th} vertex.  Here $\phi|_{\delta_i}$ denotes the restricted function where $\phi|_{\delta_i}=\phi$ at points in $\delta_i$ and vanishes elsewhere.
\label{theorem_full_decomp}
\end{theorem}

As in Theorem \ref{mainth}, if the decomposition is complete we can move from decompositions on the degree polygon level to chain level and decomposition of generic ordinarity.

\begin{cor}
Suppose the decomposition above is complete.  Then if each $f_{\Delta_i}$ is generically non-degenerate and ordinary for some prime $p$, $f$ is also generically non-degenerate and ordinary for the same prime $p$.
\end{cor}

It is often useful to use multiple decompositions simultaneously as in the parallel hyperplane decomposition.  This is possible primarily because the sum of two concave functions is itself concave.  Therefore the ``sum" of two regular decompositions is also regular.  Alternatively one may prove this more directly on the degree polygon level.  

\section{Application to Deligne Polytopes}
\label{chapter_deligne}
Using the decomposition methods established in the previous section, we are now able to investigate specific instances of generic ordinarity for a polytope $\Delta$.
We will demonstrate these methods in the case of Deligne polytopes.


\subsection{Deligne Polynomials and Polytopes}

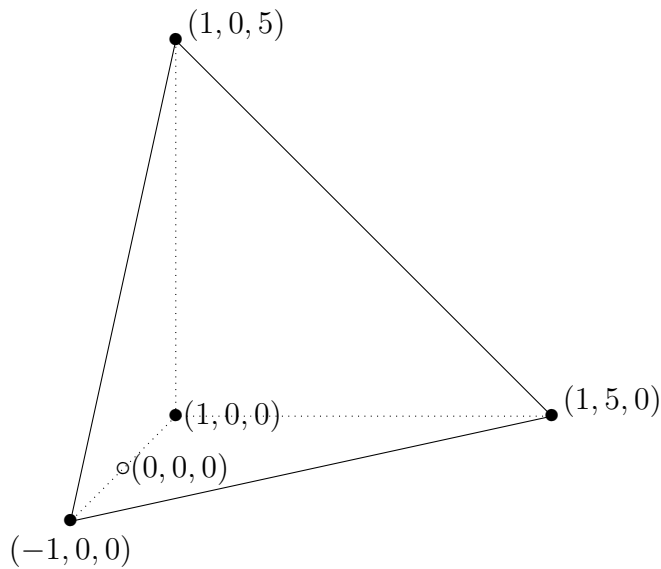
\begin{figure}[t]
\begin{center}
\begin{tikzpicture}
\draw(0,0) node {$\bullet$}; 
\draw(-0.7,-0.7) node {$\circ$}; 
\node (n1) at (-1.4,-1.4) {$\bullet$}; 
\node (v5) at (0,5) {$\bullet$}; 
\node(v0) at (5,0) {$\bullet$}; 
\draw (-1.4,-1.4) -- (0,5) -- (5,0) -- (-1.4,-1.4);
\draw[dotted] (0,5) -- (0,0) -- (5,0);
\draw[dotted] (0,0) -- (-1.4,-1.4);
\draw (n1)+(0,-0.4) node {$(-1,0,0)$};
\draw (-0.7,-0.7)+(0.75,0) node {$(0,0,0)$};
\draw (0,0)+(0.75,0) node {$(1,0,0)$};
\draw (v5)+(0.8,0.2) node {$(1,0,5)$};
\draw (v0)+(0.8,0.2) node {$(1,5,0)$};
\end{tikzpicture}
\end{center}
\caption{$\Delta$ for $d=5, n=2$}
\label{figure_deligne1}
\end{figure}

\begin{figure}[t]
\begin{center}
\begin{tikzpicture}
\draw(0,0) node {$\bullet$}; 
\draw(-0.7,-0.7) node {$\circ$}; 
\node (n1) at (-1.4,-1.4) {$\bullet$}; 
\node (v5) at (0,5) {$\bullet$}; 
\node(v0) at (5,0) {$\bullet$}; 
\draw (-0.7,-0.7) -- (0,5) -- (5,0) -- (-0.7,-0.7);
\draw[dotted] (0,5) -- (0,0) -- (5,0);
\draw[dotted] (0,0) -- (-0.7,-0.7) ;
\draw (n1)+(0,-0.4) node {$(-1,0,0)$};
\draw (-0.7,-0.7)+(0.75,0) node {$(0,0,0)$};
\draw (0,0)+(0.75,0) node {$(1,0,0)$};
\draw (v5)+(0.8,0.2) node {$(1,0,5)$};
\draw (v0)+(0.8,0.2) node {$(1,5,0)$};
\end{tikzpicture}
\end{center}
\caption{$\Delta_d'$  for $d=5, n=2$}
\label{figure_deligne2}
\end{figure}

\begin{figure}[t]
\begin{center}
\begin{tikzpicture}
\draw(0,0) node {$\bullet$}; 
\draw(-0.7,-0.7) node {$\circ$}; 
\node (n1) at (-1.4,-1.4) {$\bullet$}; 
\node (v5) at (0,5) {$\bullet$}; 
\node(v0) at (5,0) {$\bullet$}; 
\draw (-1.4,-1.4) -- (0,5) -- (5,0) -- (-1.4,-1.4);
\draw[dotted] (0,5) -- (-0.7,-0.7) -- (5,0);
\draw[dotted]  (-1.4,-1.4) -- (-0.7,-0.7) ;
\draw (n1)+(0,-0.4) node {$(-1,0,0)$};
\draw (-0.7,-0.7)+(0.75,0) node {$(0,0,0)$};
\draw (0,0)+(0.75,0) node {$(1,0,0)$};
\draw (v5)+(0.8,0.2) node {$(1,0,5)$};
\draw (v0)+(0.8,0.2) node {$(1,5,0)$};
\end{tikzpicture}
\end{center}
\caption{$\Delta_d$ for $d=5, n=2$}
\label{figure_deligne3}
\end{figure}

Let $f$ be a polynomial over $\ff_q[x_1,\ldots,x_n]$ with degree $d$ prime to $p$ and its highest degree term, say $f_d$, is a homogeneous form of degree $d$ in $n$ variables which is nonzero, and whose vanishing, if $n\geq2$, defines a smooth hypersurface in the projective space $\pp^{n-1}$.  In \cite{KATZppt}, Katz refers to such $f$ as Deligne polynomials. Following results from Browning and Heath-Brown \cite{BrHBppt}, he then examines polynomials of the form 
\begin{equation}
x_0 f(x)+g(x)+1/x_0
\label{eq_deligne}
\end{equation}
where $g$ is an arbitrary polynomial over $\ff_q$ in $n$ variables of degree $e<d/2$.  If we fix $x_0 \in \ff_q^*$ this is still a Deligne polynomial.  Interpreting $x_0$ as a variable, this is a polynomial in $n+1$ variables.  Katz gives a sharp complex estimate of the underlying exponential sum using $l$-adic cohomology.  We now consider a polytope induced by these Deligne polynomials.

Let $e_0,e_1,\ldots, e_n$ be the standard basis for $\rr^{n+1}$.  Let $\Delta$ be the polytope spanned by the origin and the vectors $-e_0, e_0,e_0+de_1+\ldots + e_0+de_n$ (see Figure \ref{figure_deligne1}).  Therefore $\Delta(x_0f+g(x)+1/x_0) =\Delta$ for generic $f$ as in (\ref{eq_deligne}).  We would like to determine conditions for generic ordinarity on $\Delta$.  There are two faces of $\Delta$ that do not contain the origin: the face $\delta_d$ spanned by $-e_0, e_0+de_1, \ldots, e_0+de_n$ and the face $\delta_d'$ spanned by $e_0,e_0+de_1, \ldots, e_0+de_n$.  The facial decomposition allows us to examine ordinarity on the individual faces $\delta_d$ and $\delta_d'$ (see Figures \ref{figure_deligne2} and \ref{figure_deligne3}).  Let $\Delta_d$ be the convex hull of $\delta_d$ and the origin.  Similarly let $\Delta_d'$ be the convex hull spanned by $\delta_d'$ and the origin.  Hence $\Delta= \Delta_d \cup \Delta_d'$.  We first examine ordinarity on $\Delta_d'$.

The face $\delta_d'$ of $\Delta_d'$ rests on the hyperplane $x_0=1$.  Therefore $D^*(\Delta_d')=1$.  Using a series of parallel hyperplanes perpendicular to the standard basis elements, Theorem 7.5 in \cite{WAN93} states that Conjecture \ref{conjAS} holds in this case.  In other words, $\Delta_d'$ is generically ordinary for every prime $p$ not dividing $d$.  

We now turn our attention to $\Delta_d$.  This behavior is far more complex.  For this reason we will call $\Delta_d$ Deligne polytopes.  The vertices of $\delta_d$ can be organized as columns in the $(n+1)\times(n+1)$-matrix:
$$\left(\begin{matrix}
-1& 1 &  1  &\ldots & 1\\
0 & d & 0 &\ldots & 0\\
0 & 0 & d & \ldots & 0\\
\vdots & & & & \vdots \\
0 & 0 & 0 &\ldots&d\end{matrix}\right).$$
Let  $e_0, e_1,\ldots, e_n$ be the standard basis in $\rr^{n+1}$.  We can write the vertices of $\delta_d$ as $-e_0,e_0+de_1,\ldots, e_0+de_n$.  The codimension 1 hyperplane spanned by $\delta_d$ is defined by the normal vector 
$$V_h=\left(\begin{matrix}
-1 \\ \frac{2}{d} \\ \vdots \\ \frac{2}{d}
\end{matrix}\right).$$
That is, for the standard inner product $\left< \star, \star \right>$, we have $\left< V_h, -e_0\right>=1$ and $\left<V_h, e_0+de_i \right>=1$ for every $1\leq i \leq h$.
Therefore 
$$D(\Delta_d)=
\begin{cases}
d & \mathrm{if~} d \mathrm{~is~odd,}\\
\frac{d}{2} & \mathrm{if~}d \mathrm{~is~even}.
\end{cases}$$  

We would like to know if $\Delta_d$ satisfies Conjecture \ref{conjAS}.  

\subsubsection{Non-Ordinarity Conditions for $\Delta_d$}
To show non-ordinarity for a fixed prime $p$, we need only to find a member $\sigma$ of the boundary decomposition that is not chain-level ordinary for any $f \in \MM_p(\Delta_d)$.  By applying the boundary decomposition theorem we may deduce that $\MM_p(\Delta_d)$ is not generically ordinary.

Consider the open line segment $\sigma$ which joins the vertices
$$
V_0= \left(\begin{matrix}
-1 \\ 0 \\ 0 \\ \vdots \\ 0
\end{matrix}\right) \mathrm{~and~}
V_1= \left(\begin{matrix}
1 \\ d \\ 0 \\ \vdots \\ 0
\end{matrix}\right).$$
The set $S(\sigma)$ as defined in Section \ref{chapter_intro} is an additive, cyclic abelian group of order $d$.  We can write a generator explicitly:
$$V_g= \left(\begin{matrix}
0 \\ 1 \\ 0 \\ \vdots \\ 0
\end{matrix}\right).$$
It follows that $w(V_g) = 1/d$.  For a prime $p$ we have $w(p V_g) = p/d \modulo{1}$.  Hence the weight function is only stable under $p$-action when $p \equiv 1 \modulo{d}$.  

Suppose $p \not\equiv 1 \modulo{d}$ and $p$ is odd.  For any $f \in \MM_p(\Delta_d)$, $f_{\sigma} = a_1/x_1 + a_2 x_1 x_2^d$ and $a_1 a_2 \neq 0$.  By Theorem \ref{theoremDiag} and Theorem \ref{th5.1} we know that the $NP(f)$ lies above $HP(\Delta_d)$.  This argument does not work for even $d$ because in that case $f_{\sigma}$ is not diagonal.

The case where $d$ is even is much more complicated, though the argument we detail below is independent of the parity of $d$.  Let $\Sigma = C(\sigma)$, the open cone generated by $\sigma$.  Recall in Section \ref{chapter_dwork}, we defined the matrix $A_1(\Sigma,f)=(a_{s,r}(f))$ with $s$ and $r$ running through $\Sigma$.  From (\ref{eq45}) and (\ref{eq49}) the coefficients of $A_1(\Sigma,f)$ are 
$$a_{s,r}(f) = F_{ps-r}(f) = F_r(f) = \sum_u \left( \prod_{j=1}^J \lambda_{u_j} a_j^{u_j}  \right) \pi^{u_1+\ldots+u_J}\pi^{w(r)-w(s)},$$
where the outer sum is over all solutions to the linear system
$$\sum_{j=1}^J u_j V_j =ps-r , u_j \geq 0, u_j \mathrm{~integral}.$$
Let $r=s=V_g$.  Then $w(ps-r)=(p-1)/D^*(\Delta_d)$.  This is an integer precisely when $p \equiv 1 \modulo {D^*(\Delta_d)}$.  From (\ref{eq46b}) we know that $(p-1)\ord F_{(p-1)V_g}(f)$ is an integer.  From (\ref{eq47}) we know the lower bound of $(p-1)\ord F_{(p-1)V_g}(f)$ is $w((p-1)V_g)$.  Therefore if $p \not\equiv 1 \modulo{D^*(\Delta_d)}$ we have
\begin{equation}
\ord F_{ps-r}(f)> w(ps-r).
\label{eq_deligneup}
\end{equation}
We can exploit (\ref{eq_deligneup}) to prove the following:
\begin{prop}
If $p \not\equiv 1 \modulo{D^*(\Delta_d)}$ then $GNP(\Delta_d,f)$ lies strictly above $HP(\Delta_d)$.  
\label{prop_nonordinary}
\end{prop}
\it Proof. \rm
Suppose $p \not\equiv 1 \modulo{D^*(\Delta_d)}$ and take $\sigma$ and $\Sigma$ as above.  The vector $V_g$ is the unique vector in $\Sigma$ with weight $1/D^*(\Delta_d)$.  Therefore $W(\Sigma,1) = 1$.  Combining this with (\ref{eq_deligneup}) we get
\begin{equation}
\ord  A_{11} >  0.
\label{eq_deligneup2}
\end{equation}
Therefore we may write $A_{11} = \xi^{\epsilon} A_{11}'$ where $\epsilon>0$ and $\ord A_{11}' \geq 0$.
We will use (\ref{eq_deligneup2}) to show that $\det(I-TA_1(\Sigma,f))$ does not coincide with its lower bound $P(\Sigma)$.
By (\ref{block53}), $A_1(\Sigma,f)$ has the block form
\begin{equation}
A_1(\Sigma,f) = \left( \begin{array}{ccccc}
A_{00} & \xi^1 A_{01} & \ldots & \xi^i A_{0i} & \ldots \\
A_{10} & \xi^1 A_{11} & \ldots & \xi^iA_{1i} & \ldots \\
\vdots & \vdots & \ddots & \vdots &\\
A_{i0} & \xi^1 A_{i1} & \ldots & \xi^iA_{ii} & \ldots \\
\vdots & \vdots & \ddots & \vdots &
\end{array} \right),
\end{equation}
where the block $A_{ij}$ is a finite matrix of $W(\Sigma,i)$ rows and $W(\Sigma,j)$ columns.  The matrix $A_{11}$ consists of a single entry. For each $i \geq 1$ one can also show that $A_{i0}$ is the $W(\Sigma,i) \times 1$ zero matrix.  Therefore 
$$\det(I-TA_1(\Sigma,f)) =
(1-TA_{00}) \det  \left( \begin{array}{ccccc}
 1-T \xi^{1+\epsilon} A_{11}' & \ldots & -T\xi^iA_{1i} & \ldots \\
 \vdots & \ddots & \vdots &\\
 -T\xi^1 A_{i1} & \ldots & I-T\xi^iA_{ii} & \ldots \\
 \vdots & \ddots & \vdots &
\end{array} \right).$$
From this we compute that the first vertex of the Newton polygon of $\det(I-TA_1(\Sigma,f))$ is 
$$(2,\frac{1}{D^*(\Delta_d)} +\epsilon).$$
This is strictly greater than the first vertex of $P(\Sigma)$ given by 
$$(2,\frac{1}{D^*(\Delta_d)}).$$
Hence we may conclude 
$$NP(\det(I-TA_1(\Sigma,f)) \gneq P(\Sigma).$$
Since $\Sigma$ is a member of the boundary decomposition we know that $NP(f) \gneq HP(\Delta_d)$ for any $f \in \MM_p(\Delta_d)$.  The proof is complete.

Note that to use this strategy to show ordinarity we would have to apply Theorem \ref{theoremDiag} to every member of the boundary decomposition.  Though the boundary decomposition is unique, this can be quite difficult to determine in high dimension.  Instead we will use decomposition theory to examine generic ordinarity.

\subsubsection{Ordinarity Conditions for $\Delta_d$}

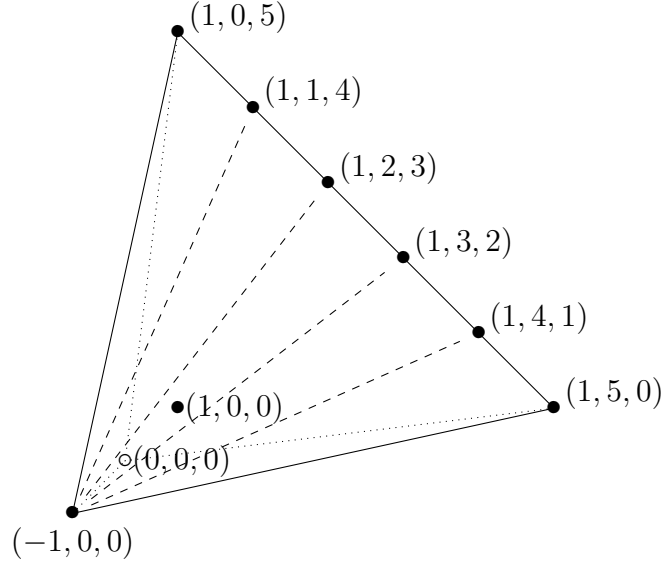
\begin{figure}[t]
\begin{center}
\begin{tikzpicture}
\draw(0,0) node {$\bullet$}; 
\draw(-0.7,-0.7) node {$\circ$}; 
\node (n1) at (-1.4,-1.4) {$\bullet$}; 
\node (v5) at (0,5) {$\bullet$}; 
\node (v4) at (1,4) {$\bullet$}; 
\node(v3) at (2,3) {$\bullet$}; 
\node(v2) at (3,2) {$\bullet$}; 
\node(v1) at (4,1) {$\bullet$}; 
\node(v0) at (5,0) {$\bullet$}; 
\draw (-1.4,-1.4) -- (0,5) -- (5,0) -- (-1.4,-1.4);
\draw[dotted] (0,5) -- (-0.7,-0.7) -- (5,0);
\draw[dotted]  (-1.4,-1.4) -- (-0.7,-0.7) ;
\draw (n1)+(0,-0.4) node {$(-1,0,0)$};
\draw (-0.7,-0.7)+(0.75,0) node {$(0,0,0)$};
\draw (0,0)+(0.75,0) node {$(1,0,0)$};
\draw (v5)+(0.8,0.2) node {$(1,0,5)$};
\draw (v4)+(0.8,0.2) node {$(1,1,4)$};
\draw (v3)+(0.8,0.2) node {$(1,2,3)$};
\draw (v2)+(0.8,0.2) node {$(1,3,2)$};
\draw (v1)+(0.8,0.2) node {$(1,4,1)$};
\draw (v0)+(0.8,0.2) node {$(1,5,0)$};
\draw[dashed] (n1) -- (v1);
\draw[dashed] (n1) -- (v2);
\draw[dashed] (n1) -- (v3);
\draw[dashed] (n1) -- (v4);
\end{tikzpicture}
\end{center}
\caption{The Decomposed Deligne Polytope for $d=5, n=2$}
\label{figure_delignedecomposed}
\end{figure}
Let $\sigma_d$ be the codimension $2$ polytope formed by the intersection of $\delta_d$ and the hyperplane defined by the points satisfying $x_0=1$:
$$\sigma_d = \delta_d \cap \{x_0=1\}.$$  Let $\sigma_d'$ be the translation of $\sigma_d$  by $-e_0$.  We can identify $\sigma_d'$ with the polytope generated by $d$ times the standard basis in $\rr^n$, in other words we are projecting $\sigma_d$ onto the last $n$ coordinates.
Any decomposition of $\sigma_d'$ will be a decomposition of $\sigma_d$.  This will in turn induce a decomposition of $\delta_d$ by taking the convex closure of $-e_0$ with each codimension $2$ sub-polytope of the decomposition of $\sigma_d$.

Consider the parallel hyperplanes $H_{i,j}$ defined by $x_i=j$ for $1\leq i \leq n$ and $0 \leq j \leq d$ defined in $\rr^n$.  This forms a decomposition $T'$ of $\sigma_d'$ into simplices, each of volume $1/n!$.  Each $H_{i,j}$ induces a hyperplane $H_{i,j}'$ in $\rr^{n+1}$ that passes through $-e_0$ and cuts $\sigma_d$.  The collection of $H_{i,j}'$ forms a triangulation of $\delta_d$.  This forms a regular decomposition $T$ of $\Delta_d$.  Since each simplex of $T'$ has volume $1/n!$ and has nonzero coordinates that sum to $d$, it follows that each member $\Delta_i$ of $T$ has $(n+1)!\Vol(\Delta_i)=d$.  

Figure \ref{figure_delignedecomposed} shows $\Delta_d$ decomposed for $d=5, n=2$.  The line segment joining $(1,0,5)$ and $(1,5,0)$ is $\sigma_d'$.  Using hyperplanes, we may decompose $\sigma_d'$ into five sub-polytopes.  These decompsitions can then be extended to $\delta_d$ by connecting each sub-polytope of $\sigma_d'$ to $(-1,0,0)$.  Finally the decomposition is extended to $\Delta_d$.  Since the decomposition of $\sigma_d'$ is a parallel hyperplane decomposition, it is regular.  Therefore the extension of this decomposition to $\delta_d$ will also be regular.

If $d$ is odd, then the intersection of $\delta_d$ and the hyperplane $\{x_1=0\}$ has no integer solutions.  Therefore $T$ is complete.  Using Theorem \ref{theoremDiag} we have
\begin{prop}
If $d$ is odd and $p \equiv 1 \modulo{d}$ then $GNP(\Delta_d,p) = HP(\Delta_d)$.  In this case Conjecture \ref{conjAS} holds.
\end{prop}
Combining this with Proposition \ref{prop_nonordinary} we have
\begin{theorem}
For odd $d$, $GNP(\Delta_d,p) = HP(\Delta_d)$ if and only if $p \equiv 1 \modulo{d}$.
\end{theorem}
We can reformulate this to apply more directly to the polynomials that Katz studies in \cite{KATZppt}.
\begin{cor}
For a fixed odd $d \in \zz_+$ and $g(x) \in \ff_q[x_1, \ldots ,x_n]$ with $\deg(g) < d/2$, let $\FF(g)$ denote the family of Laurent polynomials $tf(x)+g(x)+1/t$ parameterized by $f(x)$ of degree $d$.  Note that $\FF(g) \subset \MM_p(\Delta_d)$ for a fixed $\Delta_d$.  Conjecture \ref{conjAS} holds for this family.  Specificially if $p \equiv 1 \modulo{d}$ then,
$$\inf_{h \in \FF(g)} NP(h)  = HP(\Delta).$$
\end{cor}

If $d$ is even, then the intesection of $\delta_d$ and the hyperplane $\{x_1=0\}$ has integer solutions.  Therefore $T$ is not complete.  Subsequent decompositions are required to complete the classification of $\Delta_d$.  However direct manipulation in low dimension suggests that Conjecture \ref{conjAS} holds in this case as well.
\begin{conjecture}
$GNP(\Delta_d,p) = HP(\Delta_d)$ if and only if $p \equiv 1 \modulo{D^*(\Delta_d)}$.
\end{conjecture}
\subsection{Explicit Computation of the Hodge Polygon}

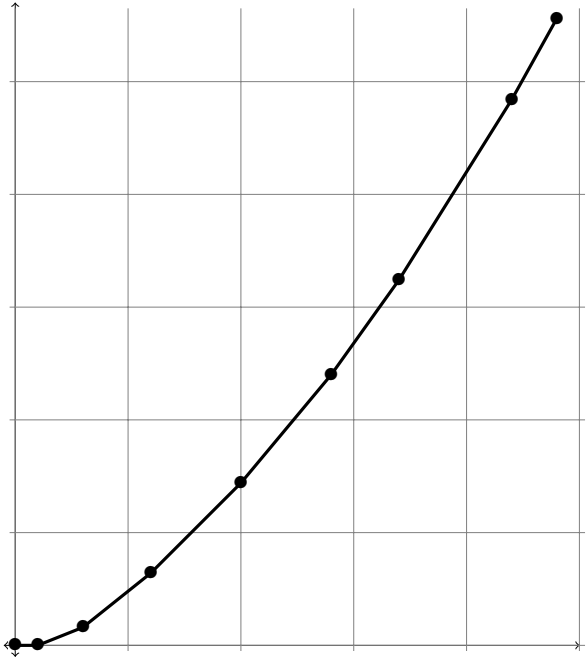
\begin{figure}[t]
\begin{center}
\begin{tikzpicture}[scale = 0.3]
\draw[to-to,very thin] (-0.5,0) -- (25,0); 
\draw[to-to,very thin] (0,-0.5) -- (0,28.5); 
\draw[step= 5cm,gray,very thin] (-0.25,-0.25) grid (25.25,28.25);
\draw[very thick] (0,0)-- (1,0)--  (3,4/5)--(6, 16/5) --(10,36/5) --(14,12)--(17,81/5) --(22,121/5) -- (24,139/5) ;
\draw (0,0) node {$\bullet$};
\draw (1,0) node {$\bullet$};
\draw (3,4/5) node {$\bullet$};
\draw(6, 16/5) node {$\bullet$};
\draw (10,36/5) node {$\bullet$};
\draw (14,12) node {$\bullet$};
\draw (17,81/5) node {$\bullet$};
\draw (22,121/5) node {$\bullet$};
\draw (24,139/5) node {$\bullet$};
\end{tikzpicture}
\end{center}
\caption{The Hodge Polygon for $\Delta_d$ Where $d=5,n=2$}
\label{figure_deligneexample}
\end{figure}

The use of linear programming and combinatorics allows for a fairly straightfoward construction of $HP(\Delta_d)$ from the definition.

Let $D=D(\Delta_d)$.  The equation for the hyperplane containing $\delta_d$ is given by 
$$(-1, \frac{2}{d},\ldots, \frac{2}{d})\smallpmatrix{u_0 \\ \vdots\\ u_n}=1.$$
Therefore we may write (\ref{eqW}) as follows:
$$W_{\Delta_d}(k) = \#\{u=(u_0,\ldots,u_n) \in L(\Delta_d) \mid (-1, \frac{2}{d},\ldots, \frac{2}{d})\smallpmatrix{u_0 \\ \vdots\\ u_n} = \frac{k}{D}\}.$$
For a given $k$ we see that $|u_0| \leq k/D$.  Solutions $u$ must satisfy the linear equation
$$-u_0 + \frac{2}{d} \sum_{i=1}^n u_i = \frac{k}{D}$$
Therefore
\begin{equation} 
\sum_{i=1}^n u_i = \frac{d}{2}(\frac{k}{D}+u_0)
\label{eqU}
\end{equation}
To contribute to $W_{\Delta}(k)$ the point $(u_0,\ldots,u_n)$ must be integral and satisfy \ref{eqU}.  The number of solutions to this equation has a standard solution in combinatorics.  Problems of this type are often called stars-and-bars problems.
Hence we get a formula 
\begin{equation}
W_{\Delta}(k) = \sum_{u_0 = -\floor{k/D}}^{\floor{k/D}}  {\frac{d}{2}(\frac{k}{D}+u_0)+(n-1)  \choose  n-1},
\label{eq_Wdeligne}
\end{equation}
where the binomial choice function is defined $0$ if either argument is not an integer.  Using this explicit formula for $W_{\Delta_d}(k)$ we are able to compute $H_{\Delta_d}(k)$ and $HP(\Delta_d)$.

\subsubsection{Example}

Using the development in this section we are able to compute the Newton polygon of $GNP(\Delta_d,p)$ for $d$ odd and $p \equiv 1 \modulo{d}$.  If $d=5, n=2$ then for $f \in \MM_p(\Delta_d)$, $L(f,T)$ is a polynomial of degree $24$ and $GNP(\Delta_d,p)=HP(\Delta_d)$.  Using (\ref{eq_Wdeligne}) and a computer, one determines that $NP(f)$ is bounded below by the lower convex hull of the vertices:
$$(1, 0),
(3, 4/5),
(6, 16/5),
(10, 36/5),
(14, 12),
(17, 81/5),
(22, 121/5),
(24, 139/5).$$
The $HP(\Delta_d)$ is displayed in Figure \ref{figure_deligneexample}.





\bibliographystyle{abbrv}
\bibliography{ple}

\end{document}